\documentclass[a4paper,11pt]{amsart}

\usepackage{array}
\usepackage{booktabs}
\usepackage{multicol}
\usepackage{multirow}
\usepackage{tabularx}
\usepackage{subfigure}
\usepackage{tikz}
\usepackage{cite}
\usepackage{soul}
\usepackage{systeme}
\usepackage{graphicx}
\usepackage[]{units}
\usepackage{color}
\usepackage{mathrsfs}
\usepackage{bm}
\usepackage[ansinew]{inputenc}
\usepackage[breaklinks,hidelinks]{hyperref}
\usepackage{float}
\usepackage{amsmath}
\usepackage{isomath}
\usepackage{amsthm}
\usepackage{amssymb}
\usepackage{amssymb}

\usepackage[foot]{amsaddr}
\usepackage[margin=2.7cm]{geometry}

\theoremstyle{remark}

\usepackage{hyperref}

\usepackage{environ}
\NewEnviron{myequation}{%
\begin{equation}
\scalebox{0.85}{$\BODY$}
\end{equation}}
\usepackage{environ}
\NewEnviron{myequation1}{%
\begin{equation}
\scalebox{0.7}{$\BODY$}
\end{equation}}
\usepackage{environ}
\NewEnviron{myequation2}{%
\begin{equation}
\scalebox{0.75}{$\BODY$}
\end{equation}}
\usepackage{environ}
\NewEnviron{myequation3}{%
\begin{equation}
\scalebox{0.9}{$\BODY$}
\end{equation}}

\usepackage{bm}
\usepackage{graphicx,color, amsfonts,amsmath}

\usepackage{mwe}

\DeclareSymbolFont{matha}{OML}{txmi}{m}{it}
\DeclareMathSymbol{\varv}{\mathord}{matha}{118}

\begin{document}

\title[Monolithic and partitioned RBM for FSI]{A monolithic and a partitioned Reduced Basis Method for Fluid--Structure Interaction problems}

\author{Monica Nonino\textsuperscript{1}}
\email{monica.nonino@univie.ac.at}

\author{Francesco Ballarin\textsuperscript{2}}
\email{fballarin@unicatt.it}

\author{Gianluigi Rozza\textsuperscript{3}}
\email{grozza@sissa.it}

\address{\textsuperscript{1}University of Vienna, Department of Mathematics.}
\address{\textsuperscript{2}Catholic University of the Sacred Heart, Department of Mathematics and Physics.}
\address{\textsuperscript{3}Sissa, mathematics area, mathLab, International School for Advanced Studies.}

\begin{abstract}
The aim of this work is to present an overview about the combination of the {Reduced Basis Method (RBM)} with two different approaches for Fluid--Structure Interaction (FSI) problems, namely a monolithic and a partitioned approach. We provide the details of implementation of two reduction procedures, and we then apply them to the same test case of interest. We first implement a reduction technique that is based on a monolithic procedure where we solve the fluid and the solid problems all at once. We then present another reduction technique that is based on a partitioned (or segregated) procedure:  the fluid and the solid problems are solved separately and then coupled using a fixed point strategy. The toy problem that we consider is based on the Turek--Hron benchmark test case, with a fluid Reynolds number $Re=100$.
\end{abstract}

\maketitle

\section{Introduction}
The bridging between approximation techniques and high-performance computing finds numerous fields of applications in the industry as well as in academia: it is sufficient to think about heat transfer problems, electromagnetic problems, structural mechanics problems (linear/nonlinear elasticity), fluid problems, and acoustic problems. In all of these examples, the models are described using a system of partial differential equations (PDE) that usually depends on a given number of {parameters} that describe the geometrical configuration of the physical domain over which the problem is formulated or that describe some physical quantities (e.g., the Reynolds number for a fluid or the Lam\'e constants for a solid) or some boundary conditions. For all of these models, we usually focus on a particular quantity of interest, also called an output of interest, such as the maximum temperature of a system, a pressure drop, or a channel flowrate. Unfortunately, computing such an output for each new value of the parameter is a difficult task that is expensive both in terms of time computation and in terms of computer memory, even on modern HPC systems. With these premises in mind, it is clear why the {Reduced Basis Method} \cite{HesthavenRozzaStamm,Haasdonk2017,Rozza2008229,BertagnaVeneziani,Lassila,Haasdonk,Nguyen} (RBM) comes into play and shows a wide range of advantages: the idea at the core of the method is to simulate the behavior of the solution of our system of interest for some chosen values of the parameters in the PDE. This is usually performed using some well-established discretization technique, such as the Finite Element Method (FEM); another discretization method used, for example, in the compressible framework in computational fluid dynamics is the Finite Volume Method (FVM), and another possibility is the Cut Finite Element Method (CutFEM) (see for example \cite{KaNoBaRo, burman_massing_hansbo, burman_hansbo}). Once we  compute these solutions, in an expensive offline phase, we can use them to build some other basis functions: with these new basis functions, in the inexpensive online phase, we can approximate the solution of the system for a new value of the parameter.\\
Among the numerous applications of the RBM, Fluid--Structure Interaction (FSI) problems {definitely} represent a great challenge as well as an extremely interesting topic (see for example \cite{BallarinRozza2016, BallarinRozzaMaday, BertagnaVeneziani, Colciago, LassilaQuarteroniRozza, LIEU20065730, DEPARIS2016700, NoBaRo19}, just to cite a few). Indeed, despite their instrinsic complicated nature (see~\cite{MadayCRAS, GRANDMONT1998525}), FSI problems are frequently used in everyday life: in naval engineering, they are used to study interactions between the water and the hull of a ship (\cite{LombardiParolini}); in biomedical applications, FSI problems are used to model the interaction between the blood flow and the deformable walls of a vessel (\cite{QuainiWang2017, Quaini:129770, BallarinJCP, String_model_1, String_model_2, Maday2009, Gerbeau2015}). Finally, in aeronautical engineering, FSI describes the way the air interacts with a plane or with (parts of) a shuttle; see \cite{FarhatPeterson, FarhatAvery, FARHAT201453, LIEU20065730}. 
The goal of this work is to present an extensive overview on the formulation of two model order reduction procedures that are applied to sets of snapshots obtained using two different approaches: a monolithic approach and a partitioned (segregated) approach. We present the entire formulation of the reduction procedures as well as some numerical results that were obtained using the two reduced order model techniques. The same test case is considered: the problem was inspired by the Turek--Hron benchmark test case FSI2, for which well-established results and analyses have already been presented in the literature; see for example \cite{turek-hron1, turek-hron2}. The main difference in our work is that, for the structure, we consider a linear problem (linearized strain tensor), different from what that considered by Turek and Hron.
The rest of the work is structured as follows: in Section~\ref{formalism}, we define the mathematical formalism behind coupled systems, and in particular, we introduce the Arbitrary Lagrangian Eulerian (ALE) formulation, which is used throughout the rest of the manuscript. In Section~\ref{two approaces}, we briefly introduce the two main approaches that can be adopted when dealing with FSI, namely monolithic and partitioned approaches. In Section~\ref{mono intro}, we present a monolithic reduced order model: in Section~\ref{mono time}, we define the time discretization; in Section~\ref{space discretization section}, we introduce the space discretization of the problem; and in \mbox{Section~\ref{lagrange multipliers intro}}, we discuss the imposition of coupling conditions through Lagrangian multipliers. In \mbox{Section~\ref{reduced basis monolithic}}, we introduce the supremizer enrichment technique, and in Section~\ref{reduced basis generation monolithic}, we present more in detail the implementation of the Proper Orthogonal Decomposition. In Section~\ref{online mono}, we formulate the online reduced order coupled system, and finally in \mbox{Section~\ref{results mono}}, we present some numerical results.
Section~\ref{part intro} is devoted to a partitioned algorithm: in Sections~\ref{part time} and~\ref{part space}, we introduce the time discretization and the space discretization, respectively. In Section~\ref{part pod}, we present the Proper Orthogonal Decomposition technique used in this case, and finally, in Sections~\ref{part online} and~\ref{part results}, we present the online reduced order problem and some numerical results. Finally, Section~\ref{discussion} is devoted to a discussion on the two algorithms.
\section{Fluid--Structure Interaction Problems}\label{formalism}
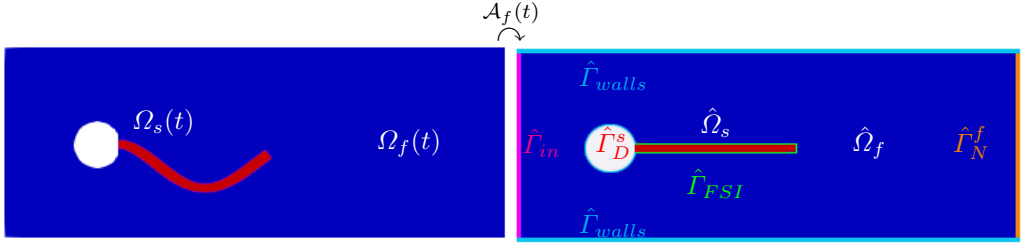
\begin{figure}
\begin{tikzpicture}
\node[anchor=south west,inner sep=0] (image) at (0,0) {\includegraphics[scale=0.55]{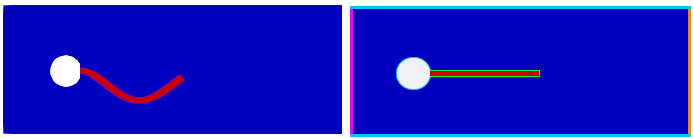}};
\begin{scope}[x={(image.south east)},y={(image.north west)}]
	\draw[white,thick] (0.16,0.6) node {\small $\Omega_s(t)$};
	\draw[white,thick] (0.4,0.5) node {\small $\Omega_f(t)$};
	\draw[black, thick] (0.5, 1.1) node {$\overset{\mathcal{A}_f(t)}{\curvearrowright}$};
	\draw[white, thick] (0.85, 0.5) node{\small$\hat{\Omega}_f$};
	\draw[white, thick] (0.7, 0.6) node{\small$\hat{\Omega}_s$};
	\draw[green, thick] (0.7, 0.3) node{\small$\hat{\Gamma}_{FSI}$};
	\draw[cyan, thick] (0.6, 0.83) node{\small$\hat{\Gamma}_{walls}$};
	\draw[cyan, thick] (0.6, 0.11) node{\small$\hat{\Gamma}_{walls}$};
	\draw[magenta, thick] (0.53, 0.5) node{\small$\hat{\Gamma}_{in}$};
	\draw[orange, thick] (0.95, 0.5) node{\small$\hat{\Gamma}_{N}^f$};
	\draw[red, thick] (0.6, 0.5) node{\small$\hat{\Gamma}_D^s$};
\end{scope}
\end{tikzpicture}
\caption{Fluid--structure interaction domains. Time-dependent configuration (\textbf{left}): fluid domain $\Omega_f(t)$ (blue) and solid deformed domain $\Omega_s(t)$ (red). Reference configuration (\textbf{right}): the inlet boundary $\hat{\Gamma}_{in}$ (magenta), the wall boundaries for the fluid $\hat{\Gamma}_{walls}$ (light blue), and the Neumann (outlet) boundary $\hat{\Gamma}_N^f$ (orange). The fluid arbitrary reference configuration $\hat{\Omega}_f$ is depicted in blue, the solid reference configuration $\hat{\Omega}_s$ is depicted in red, and $\hat{\Gamma}_D^s$ is the solid Dirichlet boundary. The fluid--structure interface in the reference configuration $\hat{\Gamma}_{FSI}$ is highlighted in green. \label{fsi_domain}}
\end{figure}

In the following, we introduce the mathematical formalism for FSI problems: we assume a two dimensional setting. Therefore, let $\Omega(t)\subset\mathbb{R}^2$ be the physical domain of interest, at time $t\in[0, T]$. The physical domain can be naturally divided into two subdomains: $\Omega(t)=\Omega_f(t)\cup\Omega_s(t)$, where $\Omega_f(t)$ is the fluid domain, and $\Omega_s(t)$, which is the solid domain; we further assume that $\Omega_f(t)\cap\Omega_s(t)=\emptyset$ and that $\bar{\Omega}_f(t)\cap\bar{\Omega}_s(t)=\Gamma_{FSI}(t)$ is the fluid--structure interface: the left side of Figure \ref{fsi_domain}  shows the fluid domain in blue and the solid domain in red. 

The fluid is assumed to be Newtonian and incompressible, and therefore, its behavior can be modelled using incompressible Navier--Stokes equation: for every $t\in[0, T]$, find $\bm{u}_f(t)\colon\Omega_f(t)\mapsto\mathbb{R}^2$ and $p_f(t)\colon\Omega_f(t)\mapsto\mathbb{R}$ such that
\begin{equation}
\label{Navier-Stokes}
\begin{cases}
\rho_f(\partial_t\bm{u}_f + ({u}_f\cdot{\nabla})\bm{u}_f)  - \text{div}\sigma_f(\bm{u}_f, p_f) = \bm{b}_f \quad\text{in }\Omega_f(t)\times(0, T],\\
\text{div}\bm{u}_f = 0 \quad \text{in }\Omega_f(t)\times(0, T],
\end{cases}
\end{equation}
where $\bm{b}_f$ is the fluid volume external force and $\sigma_f(\bm{u}_f, p_f)$ is the Cauchy stress tensor that, given the fluid is Newtonian, can be expressed in the following way:
\begin{equation*}
\sigma_f(\bm{u}_f, p_f) = \rho_f\nu_f(\nabla \bm{u}_f + \nabla^T\bm{u}_f) -p_f\boldsymbol{I}.
\end{equation*}
{where $\boldsymbol{I}$ is the $2\times2$ identity matrix, $\rho_f$ is the fluid density, and $\nu_f$ is the kinematic viscosity. We remark that, in this case, the differential operator $\nabla$ is intended to be the differentiation with respect to the time-dependent spatial coordiantes and that the same is valid for the divergence operator $div$, which also is considered with respect to the time-dependent~coordinates.}

Throughout this manuscript, for the sake of simplicity of the exposition, the solid behavior is described using a linear elasticity equation; nevertheless, we remark that everything we say can also be applied to nonlinear models for the solid. The structure problem reads as follows: for every $t\in[0, T]$, find the solid displacement $\hat{\bm{d}}_s(t)\colon\hat{\Omega}_s\mapsto\mathbb{R}^2$ such that
\begin{equation}
\label{linear elasticity}
\partial_{tt}\hat{\bm{d}}_s-\hat{\text{div}}\hat{P}(\hat{\bm{d}}_s)=\hat{\bm{b}}_s \quad\text{in }\hat{\Omega}_s\times(0, T],
\end{equation}
where $\hat{\bm{b}}_s$ is the external force acting on the solid. We assume small deformations here: this means that the solid Piola--Kirchhoff stress tensor $\hat{P}(\hat{\bm{d}}_s)$ can be expressed~as
\begin{equation*}
\hat{P}(\hat{\bm{d}}_s)= 2\mu_s\hat{\varepsilon}(\hat{\bm{d}}_s)+\lambda_s\text{tr}\hat{\varepsilon}(\hat{\bm{d}}_s)\boldsymbol{I},
\end{equation*}
with $\mu_s$ and $\lambda_s$ being the Lam\'e constants of the material and with $\hat{\varepsilon}(\hat{\bm{d}}_s)$ being the linearized strain operator, which is defined as follows:
\begin{equation*}
\hat{\varepsilon}(\hat{\bm{d}}_s):= \frac{1}{2}(\hat{\nabla}\hat{\bm{d}}_s + \hat{\nabla}^T\hat{\bm{d}}_s).
\end{equation*}

{In the previous equations, $\hat{\nabla}$ is the gradient with respect to the coordinates in the reference configuration and $\hat{div}$ is the divergence operator, always in the reference configuration.}

The first thing we notice from Equations~\eqref{Navier-Stokes} and \eqref{linear elasticity} is that the two problems are formulated over two ``different kind'' of domains: indeed, the Navier--Stokes equation is formulated over a {time-dependent}  domain $\Omega_f(t)$, whereas the linear elasticity equation is formulated over a {{time-independent}} domain $\hat{\Omega}_s$. In this section, in order to clearly introduce the Arbitrary Lagrangian Eulerian formulation, we make use of the following notation: all fields that are defined on the reference, the time-independent configuration, are denoted with a hat~$\hat{}$; on the contrary, all of the fields that are defined on a time-dependent configuration are denoted without the hat.
This distinction between the time-dependent configuration and the time-independent configuration is a rather natural peculiarity of FSI problems that arises from the different natures of the two problems involved: we refer the reader to Figure \ref{fsi_domain} for a graphical representation of both configurations.
Solving a FSI problem in a time-dependent domain can be extremely costly from a computational point of view, as it requires remeshing at every time-step in order to update the entire configuration. One possible alternative to avoid remeshing is to try instead to solve the problems in a reference configuration. For the solid problem, the definition of the reference configuration is quite easy and natural: indeed, we can say that the solid reference configuration $\hat{\Omega}_s$ is exactly the solid domain, undeformed: $\hat{\Omega}_s = \Omega_s(t=0)$. For the fluid problem, defining a reference configuration is instead rather complicated. To this aim, we introduce the ALE map
\begin{align}
\label{ALE map}
\mathcal{A}_f(t)\colon\hat{\Omega}_f&\mapsto\Omega_f(t)\nonumber\\
\hat{\bm{x}}&\mapsto \bm{x} =\hat{\bm{x}} + \hat{\bm{d}}_f(t),
\end{align}
which maps an {{arbitrary}} time-independent fluid domain to the fluid current configuration. In Equation~\eqref{ALE map}, the {{mesh displacement}} $\hat{\bm{d}}_f(t)$ at time $t$ is defined as an extension to the whole domain $\hat{\Omega}_f$ of the solid deformation $\hat{\bm{d}}_s(t)$. This extension can be carried out in different ways. Here, we adopt an harmonic extension:
\begin{equation*}
\begin{cases}
-\hat{div}(\frac{1}{J}\hat{\nabla}\hat{\bm{d}}_f) = 0, \quad \text{in }\hat{\Omega}_f\times[0, T]\\
\hat{\bm{d}}_f = \hat{\bm{d}}_s, \quad\text{on }\hat{\Gamma}_{FSI}\times[0, T].
\end{cases}
\end{equation*}
where $J$ is defined as $J(t):=\text{det}F(t)$, where $F(t)$ is the gradient of the ALE map $\mathcal{A}_f(t)$. From now on, for ease of notation, we drop the time dependence of $J$ and $F$. We complete the previous system with homogeneous Dirichlet boundary conditions on the whole $\partial\hat{\Omega}_f$.
In the previous system, we used the scaling factor $\frac{1}{J}$: the reason for doing this is that the mesh displacement obtained is sligthly more regular than the one we would obtain without the scaling factor, meaning that the deformation of the triangles of the mesh is more regular. We refer the reader to Chapter $5$ of \cite{Richter} for some interesting comparison results regarding the quality of deformation of the triangles of the mesh as well as for some
alternative ways to define the mesh displacement (linear elasticity and biharmonic extension).

Thanks to the introduction of the ALE map, we can decide to take the domain at time $t=0$ as the arbitrary time-independent fluid domain; hence, $\hat{\Omega}_f:=\Omega_f(t=0)$. After some computations (see for example \cite{Richter}), we can reformulate the Navier--Stokes problem in the reference configuration: for every $t\in [0, T]$, find the fluid velocity $\hat{\bm{u}}_f(t)\colon\hat{\Omega}_f\mapsto\mathbb{R}^2$, the fluid pressure $\hat{p}_f(t)\colon\hat{\Omega}_f\mapsto\mathbb{R}$, and the fluid displacement $\hat{\bm{d}}_f(t)\colon\hat{\Omega}_f\mapsto\mathbb{R}^2$ such that
\begin{equation}
\label{Navier-Stokes reference configuration}
\begin{cases}
\rho_fJ(\partial_t\hat{\bm{u}}_f + \hat{\nabla}\hat{\bm{u}}_fF^{-1}(\hat{\bm{u}}_f-\partial_t\hat{\bm{d}}_f))  - \hat{\text{div}}(J\hat{\sigma}_f(\hat{\bm{u}}_f, \hat{p}_f)F^{-T}) = J\hat{\bm{b}}_f \quad\text{in }\hat{\Omega}_f\times(0, T],\\
\hat{\text{div}}(J{F}^{-1}\hat{\bm{u}}_f) = 0 \quad \text{in }\hat{\Omega}_f\times(0, T],\\
-\hat{\text{div}}(\frac{1}{J}\hat{\nabla}\hat{\bm{d}}_f) = 0 \quad\text{in }\hat{\Omega}_f\times(0, T],\\
\end{cases}
\end{equation}
where $\hat{b}_f$ is now the external force acting on the fluid, expressed with respect to the coordinates in the reference configuration, and $\hat{\sigma}_f(\hat{\bm{u}}_f, \hat{p}_f)$ is the fluid Cauchy stress tensor in the arbitrary reference configuration $\hat{\Omega}_f$:
\begin{equation*}
\hat{\sigma}_f(\hat{\bm{u}}_f, \hat{p}_f):= \rho_f\nu_f(\hat{\nabla} \hat{\bm{u}}_f F^{-1} + F^{-T}\hat{\nabla}^T\hat{\bm{u}}_f).
\end{equation*}

Together, Equations~\eqref{linear elasticity} and \eqref{Navier-Stokes reference configuration} form the coupled FSI problem, expressed in a time-independent configuration. In order to have a coupled system, we still need some coupling conditions that describe the interaction between the two physics:
\begin{equation}
\label{coupling conditions}
\begin{cases}
\hat{\bm{d}}_f= \hat{\bm{d}}_s \quad \text{on }\hat{\Gamma}_{FSI}\times(0, T],\\
\hat{\bm{u}}_f = \partial_t\hat{\bm{d}}_s\quad\text{on }\hat{\Gamma}_{FSI}\times(0, T],\\
J\hat{\sigma}_fF^{-T}\hat{\bm{n}}_f = \hat{P}_s\hat{\bm{n}}_s \quad\text{on }\hat{\Gamma}_{FSI}\times(0, T],\\
\end{cases}
\end{equation}
where $\hat{\bm{n}}_f$ is the normal to the fluid--structure interface (reference configuration) from the fluid domain and $\hat{\bm{n}}_s$ is the normal to the fluid--structure interface from the solid~domain. \begin{equation*}
\begin{cases}
\hat{\bm{u}}_f = \overline{\bm{u}}(t)\quad\text{on }\hat{\Gamma}_D^f\times(0, T],\\
\hat{\bm{d}}_s = 0\quad\text{on }\hat{\Gamma}_D^s\times(0, T],\\
J\hat{\sigma}_f(\hat{\bm{u}}_f, \hat{p}_f)F^{-T}\hat{\bm{n}} = 0\quad\text{on }\hat{\Gamma}_N^f\times(0, T].
\end{cases}
\end{equation*}
where $\hat{\Gamma}_D^f:=\hat{\Gamma}_{in}\cup\hat{\Gamma}_{walls}$ is the Dirichlet boundary and is made by the inlet boundary and by $\hat{\Gamma}_{walls}$, which is the union of the top boundary, the bottom boundary, and the boundary of the cylinder immersed in the fluid. $\hat{\Gamma}_N^f$ is the fluid Neumann boundary (outflow boundary), and $\hat{\bm{n}}$ is the normal to $\hat{\Gamma}_N^f$ from the fluid domain; in addition $\overline{\bm{u}}(t)$ is the Dirichlet data, and we have 
\begin{equation}
\label{boundary velocity condition}
\overline{\bm{u}}(t)=
\begin{cases}
\bm{u}_{in}(t) \text{ on }\hat{\Gamma}_{in}\times(0, T], \\
0 \text{ on }\hat{\Gamma}_{walls}\times(0, T].
\end{cases}
\end{equation}

An example of the reference configuration together with an illustration of $\hat{\Gamma}_{FSI}$, $\hat{\Gamma}_D^f$, $\hat{\Gamma}_D^s$, $\hat{\Gamma}_N^f$, $\hat{\Gamma}_{in}$, and $\hat{\Gamma}_{walls}$ is represented in Figure \ref{fsi_domain}.

\section{Approaches to Fluid--Structure Interaction Problems}\label{two approaces}
As we previously mentioned, FSI problems are characterized by the presence of two different physics interacting with one another: we have the fluid problem, represented by the Stokes or Navier--Stokes equation, and we have the structure problem, represented by a string equation, by linear elasticity, or by nonlinear elasticity. It is therefore almost natural to conclude that there are two main different routes one can take to address such problems. 
Indeed, we can decide to solve  the fluid and the solid problems separately and then take care of the coupling between the two physics: this gives rise to the so-called partitioned (segregated) approach. 
On the other hand, we can decide to solve  the two problems together, and this gives rise to a monolithic procedure instead. 
\subsection{Monolithic Approach}
In a monolithic algorithm, 
the fluid and the solid problems are solved simultaneously. These kind of algorithms are more stable, and they usually allow for the use of bigger time-steps during the discretization of the problem in time. The main drawback is that one  really relies on {available} ad hoc softwares that are capable of handling a computational fluid dynamics problem as well as a solid mechanics problem at the same time. In addition to this, in order to pursue a Galerkin discretization of the original problem, one needs to introduce Lagrange multipliers to impose the coupling conditions at the fluid--structure interface, thus increasing the number of unknowns in the coupled problem.
For the reader interested in looking into more detail about monolithic algorithms, we refer to \cite{BallarinRozza2016, Richter, GeeWall, Quaini2017, BADIA20097986, WuCai, NoBaRo19}, even though this list is by no means complete.
\subsection{Segregated Approach}
The rationale behind a partitioned approach is to deal with the two physics separately: this indeed allows us to better exploit existing simulation tools for fluid dynamics and for structural dynamics, which are well developed nowadays and are used on a daily basis in industrial applications. 
A partitioned approach has many advantages: indeed, we have the possibility to {combine} different discretization tools for the two physics (e.g., finite volumes for the fluid and finite elements for the structure), we have the possibility to refine them for one of the two physics in time, as required by the situation. Unfortunately, in this case, there are also some drawbacks, as it turns out that, under some physical and geometrical conditions, partitioned algorithms are unstable; this situation may occur when the physical domain has a slender shape or when the fluid density $\rho_f$ is close to the solid density $\rho_s$ (this is almost always the case in hemodynamics applications, where the density of the blood is quite close to the density of the walls of the vessel). This instability occurs because of the well-known {{added mass effect} 
}: the fluid acts similar to an added mass to the solid, thus changing its natural behavior.  The reader interested in the analysis behind the added mass effect and in its derivation is referred to \cite{CausinGerbeauNobile}.  
With a partitioned procedure, we can give rise to a variety of different algorithms according to the strategy used to impose the coupling conditions at the fluid--structure interface:
\begin{itemize}
\item {{{Explicit algorithms} 
}}: after time discretization, the coupling conditions are treated explicitly at every time-step. These algorithms, also known as \emph{weakly} or \emph{loosely} coupled algorithms \cite{BUKAC2013515}, are {successfully} applied in aerodynamics applications (see \cite{FARHAT20061973, Piperno2000}), but some studies (see \cite{CausinGerbeauNobile, Gerbeau2003, LETALLEC20013039}) showed that they are unstable under some physical and geometrical conditions due to the added mass effect, as we previously mentioned. 
\item {{{Implicit algorithms}}}: in these algorithms, also known as \emph{strongly} coupled algorithms, the coupling conditions are treated implicitly at every time-step; see for example \cite{WANG20103817, Wang2018}. This implicit coupling represents a way to circumvent the instability problems due to the added mass effect; nevertheless, an implicit treatment of the coupling conditions leads to algorithms that are more expensive in terms of computational time.
\item {{{Semi-implicit algorithms}}}: in these algorithms (see \cite{BADIA20084216, Quaini2008, BASTING2017312}), the continuity of the displacement is treated explicitly whereas the other coupling conditions are treated implicitly. This alternative represents a tradeoff between the computational cost of the algorithm and its stability in relation to the physical and geometrical properties of the problem. In Section~\ref{part online}, we present a reduced order method that is based on this kind of partitioned approach.
\end{itemize}

\section{Monolithic Approach}\label{mono intro}
In the following, we propose a monolithic approach for coupled problems. As already mentioned in the Introduction, a monolithic approach means that we solve the fluid and the solid problem all together; as we  see in Section~\ref{lagrange multipliers intro}, this leads to a more subtle treatment of the coupling conditions at the fluid--structure interface; at the same time, adopting a monolithic procedure allows us to have better control of the global behavior of the coupled system.
Before going any further, we now remark that, from now on, unless otherwise stated, we assume that everything is formulated over the reference configuration. Therefore, for ease of notation, we drop the hat symbol over the variables.
\subsection{Time Discretization}\label{mono time}
We discretize the time interval $[0, T]$ with equispaced time sample points, thus obtaining $\{0=t_0,\dots, t_{N_T}=T\}$, where $t_i = i\Delta T$ and where $\Delta T$ is the time-step chosen. In the following, we  make use of the notation $f^i:=f(t^i)$ for any given function $f$.

We discretize the time derivatives in the fluid problem with a second-order backward difference formula (BDF2) (see~\cite{WuCai}):
\begin{equation*}
D_t\bm{u}_f^{i+1}=\frac{3}{2\Delta T}\bm{u}_f^{i+1}-\frac{4}{2\Delta T}\bm{u}_f^i+\frac{1}{2\Delta T}\bm{u}_f^{i-1}.
\end{equation*}

Using a Newmark scheme for the solid problem as suggested in \cite{chabannes}, we can define the structure first- and second-order time derivatives:
\begin{equation*}
\begin{split}
D_{tt}\bm{d}_s^{i+1}&=\frac{1}{\beta\Delta T}(\bm{d}_s^{i+1}-\bm{d}_s^i)-\frac{1}{\beta\Delta T}D_t\bm{d}_s^i-\Bigl(\frac{1}{2\beta}-1\Bigr)D_{tt}\bm{d}_s^i,\\
D_t\bm{d}_s^{i+1}&=\frac{\gamma}{\beta\Delta T}(\bm{d}_s^{i+1}-\bm{d}_s^i)-\Bigl(\frac{\gamma}{\beta}-1\Bigr)D_t\bm{d}_s^i-\Delta T\Bigl(\frac{\gamma}{2\beta}-1\Bigr)D_{tt}\bm{d}_s^i,
\end{split}
\end{equation*}
where the constants $\gamma$ and $\beta$ are chosen in order to ensure unconditional stability of the Newmark scheme, as suggested in \cite{chabannes}.

\subsection{Space Discretization}\label{space discretization section}
We now introduce the discrete version of the original problems \eqref{linear elasticity}--\eqref{Navier-Stokes reference configuration} from a monolithic point of view.
Let us define the following function spaces:
\begin{equation*}
\begin{split}
V^f&:=\{\bm{u}_f\in [H^1(\Omega_f)]^2\text{ s.t. }\bm{u}_f=\overline{\bm{u}}\text{ on }\Gamma^f_D\times(0, T]\},\\
V^f_0&:=\{\bm{v}_f\in [H^1(\Omega_f)]^2\text{ s.t. }\bm{v}_f=0\text{ on }\Gamma^f_D\times(0, T])\},\\
Q&:=\{q\in L^2(\Omega_f)\},\\
E^f&:=\{\bm{d}_f\in [H^1(\Omega_f)]^2\text{ s.t. }\bm{d}_f=0\text{ on }\partial\Omega_f\times(0, T])\},\\
E^s&:=\{\bm{d}_s\in [H^1(\Omega_s)]^2\text{ s.t. }\bm{d}_s=0\text{ on }\Gamma^s_D\times(0, T]\}. 
\end{split}
\end{equation*}
{We consider the spaces $V^f$, $V^f_0$, and $E^f$ to be endowed with the $H^1$ seminorm $(\nabla(\cdot), \nabla(\cdot))_{\Omega_f}$;} $Q$ to be endowed with the $L^2$ norm; and the space $E^s$ to be endowed with the $H^1$ seminorm $(\nabla(\cdot), \nabla(\cdot))_{\Omega^s}$.

We discretize the FSI problem in space, using the inf--sup stable {Taylor--Hood}  pair $(V_h^f, Q_h)=(\mathbb{P}_2, \mathbb{P}_1)$ for the fluid velocity and the fluid pressure and, similarly, for the space $V^f_0$. We use {second-order} Lagrange finite elements obtaining the discrete space $E^f_h\subset E^f$ for the mesh displacement and, similarly, for the solid displacement, resulting in the discretized space $E_h^s\subset E^s$. The weak formulation of problems \eqref{linear elasticity} and \eqref{Navier-Stokes reference configuration} now reads: find $(\bm{u}_{f, h}, {p}_{f, h}, \bm{d}_{f, h}, \bm{d}_{s, h}) \in V_h^f\times Q_h\times E_h^f\times E_h^s$ such that, for every $(\bm{v}_{f,h}, q_{f,h}, \bm{e}_{f,h}, \bm{e}_{s,h}) \in V_h^0\times Q_h\times E_h^f\times E^s_h$,
\begin{equation}
\label{coupled system}
\begin{cases}
(\partial_{tt}\bm{d}_{s,h}, \bm{e}_{s, h})_{\Omega_s} + ({P}(\bm{d}_{s,h}), \nabla \bm{e}_{s, h})_{\Omega_s}- ({P}(\bm{d}_{s,h})\bm{n}_s, \bm{e}_{s, h})_{\Gamma_{FSI}}=(\bm{b}_s, \bm{e}_{s,h})_{\Omega_s},\\
\rho_fJ(\partial_t \bm{u}_{f, h}, \bm{v}_{f, h})_{\Omega_f} + \rho_fJ(\nabla \bm{u}_{f, h}F^{-1}\bm{u}_{f, h}, \bm{v}_{f, h})_{\Omega_f}-\rho_fJ(\nabla \bm{u}_{f,h}F^{-1}\partial_t\bm{d}_{f, h}, \bm{v}_{f, h})_{\Omega_f}  + \\
+(J\sigma_f(\bm{u}_{f,h}, p_{f, h})F^{-T}, \nabla \bm{v}_{f, h})_{\Omega_f} - (J\sigma_f(\bm{u}_{f,h}, p_{f, h})F^{-T}\bm{n}_f, \bm{v}_{f,h})_{\Gamma_{FSI}}= (J\bm{b}_f, \bm{v}_{f,h})_{\Omega_f}\\
-(\text{div}(JF^{-1}\bm{u}_{f, h}), q_{f,h})_{\Omega_f} = 0 \\
(\frac{1}{J}\nabla\bm{d}_{f,h}, \nabla\bm{e}_{f,h})_{\Omega_f} = (\frac{1}{J}\nabla \bm{d}_{f, h}\bm{n}_f, \bm{e}_{f, h})_{\Gamma_{FSI}}.\\
\end{cases}
\end{equation}

{In the previous system, $\bm{n}_f$ and $\bm{n}_s$ are the normals to the fluid--structure interface, from the fluid domain and the solid domain, respectively.}
\subsubsection{Coupling Conditions through Lagrange Multipliers}\label{lagrange multipliers intro}
Let us analyze the boundary integrals appearing in \mbox{system~\eqref{coupled system}} a little bit more in detail. From condition \eqref{coupling conditions}, {it} 
 is clear that we have to satisfy the following:
\begin{equation*}
(P(\bm{d}_{s,h})\bm{n}_s, \bm{e}_{s, h})_{\Gamma_{FSI}} = -(J\sigma_f(\bm{u}_{f,h}, p_{f,h})F^{-T}\bm{n}_f, \bm{v}_{f,h})_{\Gamma_{FSI}}.
\end{equation*}

Now, let us have a look at {condition} 
 \eqref{coupling conditions}. It is easy to see that, in the time-continuous regime, thanks {to} \eqref{coupling conditions}, the following are equivalent:
\begin{equation*}
\bm{d}_{f,h}=\bm{d}_{s, h} \text{ on }\Gamma_{FSI}, \quad \partial_t\bm{d}_{s,h}=\bm{u}_{f,h}\text{ on }\Gamma_{FSI}, \quad\partial_t\bm{d}_{s,h}=\partial_t\bm{d}_{f,h}\text{ on }\Gamma_{FSI}.
\end{equation*}

This equivalence does not hold anymore in general after time discretization, as one may choose different time integration schemes for the fluid and for the solid. In this work, we chose to weakly enforce the two following conditions:
\begin{equation}
\bm{d}_{f,h}=\bm{d}_{s,h} \text{ on }\Gamma_{FSI}, \quad \partial_t\bm{d}_{s,h}=\bm{u}_{f,h}\text{ on }\Gamma_{FSI}.
\end{equation}

In order to do so, let us introduce the space $L:=[H^{-\frac{1}{2}}(\Gamma_{FSI})]^2$ first and let us consider a finite dimensional subspace $L_h\subset L$ (we use first-order Lagrange finite elements). Let us take a discretized Lagrangian multiplier field $\bm{\lambda}_{u,h}\in L_h$ and its corresponding test function $\bm{\mu}_{u, h}\in L_h$. We can then write
\begin{equation*}
(\bm{\mu}_{u, h}, \bm{u}_{f,h}-\partial_t\bm{d}_{s,h})_{\Gamma_{FSI}}=0,
\end{equation*}
and, by identifying the Lagrange multiplier with the unknown surface traction, we can finally rewrite the surface integrals in \eqref{coupled system} 
 as follows:
\begin{equation*}
\begin{split}
- ({P}(\bm{d}_{s,h})\bm{n}_s, \bm{e}_{s, h})_{\Gamma_{FSI}} &= -(\bm{\lambda}_{u,h}, \bm{e}_{s,h})_{\Gamma_{FSI}}, \\
-(J\sigma_f(\bm{u}_{f,h}, p_{f, h})F^{-T}\bm{n}_f, \bm{v}_{f,h})_{\Gamma_{FSI}} &= +(\bm{\lambda}_{u,h}, \bm{e}_{s,h})_{\Gamma_{FSI}}.
\end{split}
\end{equation*} 

We treat the continuity of the displacements at the interface similarly: we take another Lagrangian multiplier field $\bm{\lambda}_{d,h}\in L_h$ and its corresponding test function $\bm{\mu}_{d, h}\in L_h$. We can then write
\begin{equation*}
(\bm{\mu}_{d, h}, \bm{d}_{f,h}-\bm{d}_{s,h})_{\Gamma_{FSI}}=0,
\end{equation*}
and now we identify the Lagrange multiplier with the surface traction caused by the mesh displacement; thus, we can rewrite the surface integral in \eqref{coupled system}$_4$ as follows:
\begin{equation*}
(\frac{1}{J}\nabla \bm{d}_{f, h}\bm{n}_f, \bm{e}_{f, h})_{\Gamma_{FSI}} = (\bm{\lambda}_{d,h}, \bm{e}_{f,h})_{\Gamma_{FSI}}.
\end{equation*} 

Finally, after the space and the time discretization, our monolithic coupled system reads as follows: for every $t^i$, $i=0, \dots, N_T$, find $\bm{u}_{f, h}^{i+1}\in V_h^f$, $p_{f, h}^{i+1}\in Q_h$, $\bm{d}_{f, h}^{i+1}\in E_h^f$, $\bm{d}_{s, h}^{i+1} \in E_h^s$, $\bm{\lambda}_{u, h}^{i+1}\in L_h$, and $\bm{\lambda}_{d, h}^{i+1}\in L_h$ such that, for every $\bm{v}_{f,h}\in V_h^0$, $q_{f,h}\in Q_h$, $\bm{e}_{f,h}\in E_h^f$, $\bm{e}_{s,h}\in E^s_h$, $\bm{\mu}_{u, h}\in L_h$, and $\bm{\mu}_{d, h}\in L_h$,
\begin{equation}
\label{coupled discretized system}
\begin{cases}
(D_{tt}\bm{d}_{s,h}^{i+1}, \bm{e}_{s, h})_{\Omega_s} + ({P}(\bm{d}_{s,h}^{i+1}, \nabla \bm{e}_{s, h})_{\Omega_s}+ (\bm{\lambda}_{u,h}^{i+1}, \bm{e}_{s, h})_{\Gamma_{FSI}}=(\bm{b}_s, \bm{e}_{s,h})_{\Omega_s},\\
\rho_fJ(D_t \bm{u}_{f, h}^{i+1},\bm{v}_{f, h})_{\Omega_f} + \rho_fJ(\nabla \bm{u}_{f, h}^{i+1}F^{-1}\bm{u}_{f, h}^{i+1}, \bm{v}_{f, h})_{\Omega_f}-\rho_fJ(\nabla \bm{u}_{f,h}F^{-1}D_t\bm{d}_{f, h}^{i+1}, \bm{v}_{f, h})_{\Omega_f}  + \\
+(J\sigma_f(\bm{u}_{f,h}^{i+1}, p_{f, h}^{i+1})F^{-T}, \nabla \bm{v}_{f, h})_{\Omega_f} - (\bm{\lambda}_{u, h}^{i+1}, \bm{v}_{f,h})_{\Gamma_{FSI}}= (J\bm{b}_f, \bm{v}_{f,h})_{\Omega_f}\\
-(\text{div}(JF^{-1}\bm{u}_{f, h}^{i+1}), q_{f,h})_{\Omega_f} = 0 \\
(\frac{1}{J}\nabla\bm{d}_{f,h}^{i+1}, \nabla\bm{e}_{f,h})_{\Omega_f} = (\bm{\lambda}_{d, h}^{i+1}, \bm{e}_{f, h})_{\Gamma_{FSI}},\\
(\bm{u}_{f, h}^{i+1}- D_t\bm{d}_{s, h}^{i+1}, \bm{\mu}_{u, h})_{\Gamma_{FSI}}=0,\\
(\bm{d}_{f, h}^{i+1} - \bm{d}_{s, h}^{i+1}, \bm{\mu}_{d, h})_{\Gamma_{FSI}}=0.
\end{cases}
\end{equation}
\subsection{Lifting Function and Supremizer Enrichment}\label{reduced basis monolithic}
{We now present the technique used to perform a compression on the set of snapshots obtained with the monolithic FE discretization in order to create a set of reduced basis functions. The first detail that we explain is the introduction of a lifting function for the fluid velocity $\bm{u}_{f, h}$. The use of a lifting function is quite common in the RBM approach; see for example \cite{HesthavenRozzaStamm, Supremizer}: the advantage of this technique is represented by the fact that, sometimes, in the problem of interest, we have to deal with non-homogeneous Dirichlet boundary conditions, as in our case, where we remind the reader that we ask for $\bm{u}_{f,h}(t)=\bm{u}_{in}(t)$ at the inlet boundary for every $t\in(0, T]$. From the implementation point of view, this condition does not present any problems during the offline phase, where we solve the FE discretization of the original problem. However, this non-homogeneous condition represents a problem during the solution of the online system. Indeed, if we perform a POD on a set of snapshots that satisfy the same inlet condition, we obtain a set of reduced basis functions for the fluid velocity that has a given value at the inlet boundary. Now, let us assume that we have solved the online system and, hence, have found the coefficients that allow us to represent the reduced velocity approximation as a linear combination of the velocity reduced basis functions. It is easy to imagine that a linear combination of these basis functions does not automatically satisfy the original inlet condition in general.
A solution to this problem is represented by the introduction of a lifting function $\bm{\ell}_h(t)\in V^f_h$  during the offline phase such that $\bm{\ell}_h(t) = \bm{u}_{in}(t)$ on $\Gamma_D^f$ for every $t\in(0, T]$. By subtracting the lifting function to the fluid velocity snapshots $\bm{u}_{f, h}$ before performing a POD, we obtain a new variable $\bm{u}_{0,h}:=\bm{u}_{f, h}-\bm{\ell}_h \in V_0^f$ that satisfies a homogeneous Dirichlet boundary condition also at the fluid inlet boundary $\Gamma_{in}$: these are the snapshots on which we perform a POD, thus obtaining basis functions that are all zero at the inlet boundary.
We point out that the definition of the lifting function is not unique; for our problem, the lifting function $\bm{\ell}_h$ is defined by solving a steady Stokes problem over $\Omega_f$ at every time-step, with a prescribed inlet velocity $\bm{u}_{in}$ and with homogeneous Dirichlet boundary conditions on $\Gamma_{walls}\cup\Gamma_{FSI}$, but other choices are possible as well.}

Another important technique that has been used in this manuscript is the supremizer enrichment technique, which is necessary in order to obtain a stable approximation of the fluid pressure also at the reduced order level. The main reason for the introduction of the supremizer enrichment is that, even if the FE spaces $(V_h^f, Q_h)$ satisfy the inf--sup condition {(which guarantees that the Navier--Stokes problem is uniquely solvable with respect to the pressure; see for example \cite{boffi_mixed,BREZZI199027})}, this may not hold true anymore once we move to the reduced spaces generated by the reduced basis functions. {With the introduction of the supremizer enrichment, therefore, we aim to construct a pair of reduced function spaces for the fluid velocity and the fluid pressure, that also satisfies the \emph{inf--sup} condition.}
The \emph{supremizer variable} $\bm{s}_h\in V_h^0$, is defined by solving the following problem: find $\bm{s}_h\in V_h^0$ such~that
\begin{equation*}
-(\text{div}\bm{v}_{f,h}, p_{f, h})_{\Omega_f} = (\nabla \bm{s}_h, \nabla \bm{v}_{f,h})_{\Omega_f} \quad \forall \bm{v}_{f, h}\in V_h^0.
\end{equation*} 

In the previous equation, $p_{f, h}$ is the FE pressure solution of the Navier--Stokes problem whereas the right-hand side is the scalar product that defines the $H^1$ seminorm, which we consider for the velocity function space $V^0_h$. For the reader interested in the details about this technique as well as the computations that lead to the formulation of the previous problem, we refer to \cite{Supremizer, shafqat}.

\subsection{Reduced Basis Generation}\label{reduced basis generation monolithic}
Once we obtain the FE supremizer snapshots $\bm{s}_h^i$, $i=0,\dots, N_T$, and once we homogenize the fluid velocity snapshots thanks to the lifting function, we are ready to generate a set of reduced basis by performing a compression by Proper Orthogonal~Decomposition.

{In order to perform a POD, we need two main ingredients: the matrices of the inner products and the snapshots matrices. First, we need to introduce the basis functions for the FE spaces that we consider. We define the following:}
\begin{eqnarray*}
\{\bm{\varphi}_1^u,\dots,\bm{\varphi}_{\mathcal{N}_h^u}^u\} \text{ the FE basis of the discretized space } V^0_h,\nonumber\\
\{{\varphi}_1^p,\dots,{\varphi}_{\mathcal{N}_h^p}^p\} \text{ the FE basis of the discretized space } Q_h,\nonumber\\
\{\bm{\varphi}_1^{d_f},\dots,\bm{\varphi}_{\mathcal{N}_h^{d_f}}^{d_f}\} \text{ the FE basis of the discretized space } E^f_h,\nonumber\\
\{\bm{\varphi}_1^{d_s},\dots,\bm{\varphi}_{\mathcal{N}_h^{d_s}}^{d_s}\} \text{ the FE basis of the discretized space } E^s_h,\nonumber\\
\{\bm{\varphi}_1^{\lambda},\dots,\bm{\varphi}_{\mathcal{N}_h^{\lambda}}^{\lambda}\} \text{ the FE basis of the discretized space } L_h,
\end{eqnarray*}
{where $\mathcal{N}_h^u$ is the dimension of the FE space $V^0_h$, $\mathcal{N}_h^p$ is the dimension of the FE space $Q_h$, $\mathcal{N}_h^{d_f}$ is the dimension of the FE space $E^f_h$, $\mathcal{N}_h^{d_s}$ is the dimension of the FE space $E^s_h$, and $\mathcal{N}_h^{\lambda}$ is the dimension of the FE space $L_h$ (which we remember we used to approximate both the Lagrange multiplier $\bm{\lambda}_u$ and the Lagrange multiplier $\bm{\lambda}_d$).
We begin by constructing the snapshots matrices $\mathcal{S}_u\in\mathbb{R}^{\mathcal{N}_h\times M}$, $\mathcal{S}_s\in\mathbb{R}^{\mathcal{N}_h\times M}$,$\mathcal{S}_p\in\mathbb{R}^{\mathcal{N}_h\times M}$, $\mathcal{S}_{d_f}\in\mathbb{R}^{\mathcal{N}_h\times M}$, $\mathcal{S}_{d_s}\in\mathbb{R}^{\mathcal{N}_h\times M}$, $\mathcal{S}_{\lambda_u}\in\mathbb{R}^{\mathcal{N}_h\times M}$, and $\mathcal{S}_{\lambda_d}\in\mathbb{R}^{\mathcal{N}_h\times M}$ defined as follows:}
\begin{eqnarray}
\mathcal{S}_u&&=[\bm{u}_{0,h}^1,\dots,\bm{u}_{0, h}^{N_T}, 0, \dots, 0, 0,\dots,0,0,\dots,0,0\dots,0,0,\dots,0], \nonumber\\
\mathcal{S}_s&&{=[\bm{s}_{h}^1,\dots,\bm{s}_{h}^{N_T}, 0, \dots, 0, 0,\dots,0,0,\dots,0,0\dots,0,0,\dots,0],} \nonumber\\
\mathcal{S}_p&&{=[0, \dots, 0, p_{f,h}^1,\dots,p_{f, h}^{N_T}, 0, \dots, 0, 0, \dots, 0, 0,\dots,0,0,\dots,0,0\dots,0], }\nonumber\\
\mathcal{S}_{d_f}&&{=[0, \dots, 0, 0, \dots, 0, \bm{d}_{f, h}^1,\dots, \bm{d}_{f,h}^{N_T}, 0,\dots, 0, 0, \dots, 0, 0,\dots,0], }\nonumber\\
\mathcal{S}_{d_s}&&{=[0, \dots, 0, 0,\dots,0,0,\dots,0,\bm{d}_{s, h}^1,\dots, \bm{d}_{s,h}^{N_T}, 0,\dots, 0,0,\dots,0],}  \nonumber\\
\mathcal{S}_{{\lambda}_u}&&{=[0, \dots, 0, 0,\dots,0, 0,\dots,0,0,\dots,0, \bm{\lambda}_{u, h}^1,\dots, \bm{\lambda}_{u,h}^{N_T}, 0,\dots, 0],}  \nonumber\\
\mathcal{S}_{{\lambda}_d}&&{=[0, \dots, 0, 0,\dots,0, 0,\dots,0,0,\dots,0,0,\dots,0\bm{\lambda}_{d, h}^1,\dots, \bm{\lambda}_{d,h}^{N_T}]. }\nonumber
\end{eqnarray}

{In the previous definition, we have that $\mathcal{N}_h= \mathcal{N}_h^u+\mathcal{N}_h^p+\mathcal{N}_h^{d_f}+\mathcal{N}_h^{d_s}+2\mathcal{N}_h^{\lambda}$ is the sum of the dimensions of all  FE spaces that we  used for the FE approximation of each component of the solution of the FSI system and that $M=6N_T$.
{Next, we need to define the inner product matrices $X_u$, $X_p$, $X_{d_s}$, $X_{d_f}$, $X_{\lambda_u}$, and $X_{\lambda_d}$, all belonging to $\mathbb{R}^{\mathcal{N}_h\times\mathcal{N}_h}$. These matrices are block diagonal matrices and have the following~form:}}
\begin{eqnarray*}
X_u = \text{diag}(\bm{x}_u, \bm{0}_p, \bm{0}_{{d}_f}, \bm{0}_{{d}_s}, \bm{0}_{{\lambda}_u}, \bm{0}_{{\lambda_d}})\nonumber \\
X_p = \text{diag}(\bm{0}_u, \bm{x}_p, \bm{0}_{{d}_f}, \bm{0}_{{d}_s}, \bm{0}_{{\lambda}_u}, \bm{0}_{{\lambda_d}})\nonumber \\
X_{d_f} =  \text{diag}(\bm{0}_u, \bm{0}_p, \bm{x}_{{d}_f}, \bm{0}_{{d}_s}, \bm{0}_{{\lambda}_u}, \bm{0}_{{\lambda_d}})\nonumber \\
X_{d_s} =  \text{diag}(\bm{0}_u, \bm{0}_p, \bm{0}_{{d}_f}, \bm{x}_{{d}_s}, \bm{0}_{{\lambda}_u}, \bm{0}_{{\lambda_d}})\nonumber \\
X_{\lambda_u} =  \text{diag}(\bm{0}_u, \bm{0}_p, \bm{0}_{{d}_f}, \bm{0}_{{d}_s}, \bm{x}_{{\lambda}_u}, \bm{0}_{{\lambda_d}})\nonumber \\
X_{\lambda_d} =  \text{diag}(\bm{0}_u, \bm{0}_p, \bm{0}_{{d}_f}, \bm{0}_{{d}_s}, \bm{0}_{{\lambda}_u}, \bm{x}_{{\lambda_d}}).\nonumber \\
\end{eqnarray*}

{In the previous definition, for simplicity, we used the following notation: $\bm{0}_*\in\mathbb{R}^{\mathcal{N}_h^*\times\mathcal{N}_h^*}$ is a zero block of dimension $\mathcal{N}_h^*\times\mathcal{N}_h^*$, where $*\in\{u_f, p, d_f, d_s, \lambda_u, \lambda_d\}$. In addition to this, we have the following nonzero blocks:}
\begin{eqnarray*}
&&{(\bm{x}_u)_{i,j}= 
(\nabla\bm{\varphi}_i^u, \nabla\bm{\varphi}_j^u)_{\Omega_f} \text{ for }i,j=1,\dots,\mathcal{N}_h^u,}
\nonumber\\
&&{(\bm{x}_p)_{i,j} = 
(\varphi_i^p, \varphi_j^p)_{\Omega_f} \text{ for }i,j=1,\dots,\mathcal{N}_h^p,}
\nonumber\\
&&{(\bm{x}_{d_f})_{i,j}= 
(\nabla\bm{\varphi}_i^{d_f}, \nabla\bm{\varphi}_j^{d_f})_{\Omega_f} \text{ for }i,j=1,\dots,\mathcal{N}_h^{d_f},}
\nonumber\\
&&{(\bm{x}_{d_s})_{i,j} = 
(\nabla\bm{\varphi}_i^{d_s}, \nabla\bm{\varphi}_j^{d_s})_{\Omega_s} \text{ for }i,j=1,\dots,\mathcal{N}_h^{d_s},}
\nonumber\\
&&{(\bm{x}_{\lambda_u})_{i,j}= 
(\bm{\varphi}_i^{\lambda}, \bm{\varphi}_j^{\lambda})_{\Gamma_{FSI}} \text{ for }i,j=1,\dots,\mathcal{N}_h^{\lambda},}
\nonumber\\
&&{(\bm{x}_{\lambda_d})_{i,j} = 
(\bm{\varphi}_i^{\lambda}, \bm{\varphi}_j^{\lambda})_{\Gamma_{FSI}} \text{ for }i,j=1,\dots,\mathcal{N}_h^{\lambda}.}
\end{eqnarray*}

{As the reader may notice, the inner product matrices are very big, given the fact that $\mathcal{N}_h^u,\dots,\mathcal{N}_h^{\lambda}\gg1$: this is because the structure of the Proper Orthogonal Decomposition itself reflects the fact that we use a monolithic approach to solve the FSI problem.

We are now able to define the correlation matrices $\mathcal{C}_u$,  $\mathcal{C}_s$,  $\mathcal{C}_p$,  $\mathcal{C}_{d_f}$,  $\mathcal{C}_{d_s}$,  $\mathcal{C}_{\lambda_u}$, and $\mathcal{C}_{\lambda_d}$, all belonging to the space $\mathbb{R}^{M\times M}$:}
\begin{eqnarray*}
\mathcal{C}_u&:=\mathcal{S}_u^TX_u\mathcal{S}_u\nonumber\\
\mathcal{C}_s&:=\mathcal{S}_s^TX_u\mathcal{S}_s\nonumber\\
\mathcal{C}_p&:=\mathcal{S}_p^TX_p\mathcal{S}_p\nonumber\\
\mathcal{C}_{d_f}&:=\mathcal{S}_{d_f}^TX_{d_f}\mathcal{S}_{d_f}\nonumber\\
\mathcal{C}_{d_s}&:=\mathcal{S}_{d_s}^TX_{d_s}\mathcal{S}_{d_s}\nonumber\\
\mathcal{C}_{\lambda_u}&:=\mathcal{S}_{\lambda_u}^TX_{\lambda_u}\mathcal{S}_{\lambda_u}\nonumber\\
\mathcal{C}_{\lambda_d}&:=\mathcal{S}_{\lambda_d}^TX_{\lambda_d}\mathcal{S}_{\lambda_d}\nonumber\\
\end{eqnarray*}

{\emph{Remark}: in the correlation matrices, all of the snapshots are defined on a common mesh. Indeed, as we mentioned at the beginning of the section, we dropped the hat notation~$\hat{}$, with the understanding that all of the quantities are defined on the common reference configuration. This aspect is extremely important, as mapping everything back onto a reference configuration greatly simplifies the implementation of the Proper Orthogonal Decomposition.}

{Once we built the correlation matrices, we carried out a POD compression on the set of snapshots, following for example \cite{Kunisch2002492}. We did this by solving the following (seven) eigenvalue problems: 
\begin{equation}
\label{eigenvalue problem}
\mathcal{C}_{*}\mathcal{Q}_{*}=\mathcal{Q}_{*}\Lambda_{*},
\end{equation}
where $*\in\{u_f$, $s$, $p$, $d_f$, $d_s$, $\lambda_u$, $\lambda_d\}$, $\mathcal{Q}_*$ is the eigenvectors matrix, and $\Lambda_*$ is the diagonal eigenvalues matrix. The $k$th reduced basis function related to problem \eqref{eigenvalue problem} is obtained by applying the snapshots matrix $\mathcal{S}_*$ to the $k$th column of the matrix $\mathcal{Q}_*$; we therefore end up with the following basis functions (assume for simplicity $*=u_f$):
\begin{equation*}
\bm{\Phi}_k^u:= \frac{1}{\sqrt{\lambda_k^u}}\mathcal{S}_u\underline{\bm{v}}_k^u,
\end{equation*}
where $\lambda_k^u$ is the eigenvalue corresponding to the eigenvector $\bm{v}_k^u$.
Similar definitions hold true for the other components of the solution, namely $p_f$, $\bm{d}_s$, $\bm{d}_f$, $\bm{\lambda}_u$, and $\bm{\lambda}_d$, as well as for the supremizer $\bm{s}$. We refer the reader interested in more details about the POD to \cite{HesthavenRozzaStamm, Supremizer}. }

{We therefore end up with the following set of reduced basis: $\{\bm{\Phi}_1^{s,u}, \dots, \bm{\Phi}_{N_u}^{s,u}, \dots, \bm{\Phi}_1^{\lambda_d}, \dots,$\linebreak$ \bm{\Phi}_{N_{\lambda_d}}^{\lambda_d}\}$}, where each basis function is a block function of six components (one for each variable in the FSI problem):
\begin{equation*}
\bm{\Phi}_k^{s, u}=
\begin{pmatrix}
\Phi_k^{s,u}\\
0\\
0\\
0\\
0\\
0
\end{pmatrix},
\dots,
\bm{\Phi}_k^{\lambda_d}=
\begin{pmatrix}
0\\
0\\
0\\
0\\
0\\
\Phi_k^{\lambda_d}
\end{pmatrix}.
\end{equation*}
{where we have used the following notation: let $\{\bm{\Phi}_1^u,\dots,\bm{\Phi}_{N1}^u\}$ be the basis functions obtained by a compression by POD on the fluid velocity snapshots and let $\{\bm{\Phi}_1^s,\dots,\bm{\Phi}_{N2}^s\}$ be the basis functions obtained by running a POD on the supremizer snapshots. In order to carry on the supremizer enrichment technique, we consider the union of the two sets of basis functions $\{\bm{\Phi}_1^u,\dots,\bm{\Phi}_{N1}^u, \bm{\Phi}_1^s,\dots,\bm{\Phi}_{N2}^s\}$ and then denote by $\bm{\Phi}_k^{s,u}$ a generic element of the last set. We indicated $\bm{\Phi}_k^{s,u}$ in order to remark that the reduced basis functions for the fluid velocity consist of the reduced basis generated by the fluid velocity snapshots and of the reduced basis generated by the supremizer snapshots. 

{Finally, we introduce the reduced order finite dimensional space $\bold{V}_N=\text{span}\{\bm{\Phi}^{s, u}_1,\dots,$\linebreak$\bm{\Phi}_{N_u}^{s,u}, \dots, \bm{\Phi}_1^{\lambda_d}, \dots, \bm{\Phi}_{N_{\lambda_d}}^{\lambda_d}\}$,} with $N_u=N1+N2$ and $N=N_u+\dots+N_{\lambda_d}$.}
\subsection{Online Phase}\label{online mono}
Once we have the reduced basis functions, we can define the reduced solution $\mathbf{u}_N^{i+1}:=(\bm{u}_{0, N_u}^{i+1}, p_{N_p}^{i+1}, \bm{d}_{{N_{d_f}}}^{i+1}, \bm{d}_{N_{d_s}}^{i+1}, \bm{\lambda}_{N_{\lambda_u}}^{i+1}, \bm{\lambda}_{N_{\lambda_d}}^{i+1})$ of our FSI problem:
\begin{align}
\bm{u}_{0, N_u}^{i+1} &:=\sum_{k=0}^{N_u}\underline{u}_{0,k}^{i+1}\Phi_k^{s, u},\\
p_{N_p}^{i+1} &:=\sum_{k=0}^{N_p}\underline{p}_k^{i+1}\Phi_k^p,\\
\bm{d}_{N_{d_f}}^{i+1} &:=\sum_{k=0}^{N_{d_f}}\underline{d}_{f, k}^{i+1}\Phi_k^{d_f},\\
\bm{d}_{N_{d_s}}^{i+1} &:=\sum_{k=0}^{N_{d_s}}\underline{d}_{s, k}^{i+1}\Phi_k^{d_s},\\
\bm{\lambda}_{N_{\lambda_u}}^{i+1} &:=\sum_{k=0}^{N_{\lambda_u}}\underline{\lambda}_{u, k}^{i+1}\Phi_k^{\lambda_u},\\
\bm{\lambda}_{N_{\lambda_d}}^{i+1} &:=\sum_{k=0}^{N_{\lambda_d}}\underline{\lambda}_{d, k}^{i+1}\Phi_k^{\lambda_d}.
\end{align}

In the previous equations, the underline bar indicates the vector of coefficients of the reduced solution; therefore, it indicates an element of $\mathbb{R}$ (for scalar components of the reduced solution) or $\mathbb{R}^2$ for vectorial components of the reduced solution (such as the fluid velocity for example).
The online monolithic reduced order system reads as follows: for every $t^{i+1}$, $i=0,\dots, N_T-1$, {find $\bold{u}_N^{i+1}\in \bold{V}_N$, $\bold{u}_N^{i+1} =(\bm{u}_{0, N_u}^{i+1}, p_{N_p}^{i+1}, \bm{d}_{N_{d_f}}^{i+1}, \bm{d}_{N_{d_s}}^{i+1}, \bm{\lambda}_{N_{\lambda_u}}^{i+1}, \bm{\lambda}_{N_{\lambda_d}}^{i+1})$ such that, for all $\bold{v}_N\in \bold{V}_N$, $\bold{v}_N=(\bm{v}_{N_u}, q_{N_p}, \bm{e}_{N_{d_f}}, \bm{e}_{N_{d_s}}, \bm{\mu}_{N_{\lambda_u}}, \bm{\mu}_{N_{\lambda_d}})$, the following  holds}
\begin{equation}
\label{online system monolithic}
\begin{cases}
(D_{tt}\bm{d}_{N_{d_s}}^{i+1}, \bm{e}_{N_{d_s}})_{\Omega_s} + ({P}(\bm{d}_{N_{d_s}}^{i+1}), \nabla\bm{e}_{N_{d_s}})_{\Omega_s} + (\bm{\lambda}_{N_{\lambda_u}}^{i+1}, \bm{e}_{N_{d_s}})_{\Gamma_{FSI}}=(\bm{b}_s, \bm{e}_{N_{d_s}})_{\Omega_s},\\
\rho_fJ(D_t(\bm{u}_{0, N_u}^{i+1}+\bm{\ell}^{i+1}_{N_u}), \bm{v}_{N_u})_{\Omega_f} + \rho_fJ(\nabla(\bm{u}_{0, N}^{i+1}+\bm{\ell}^{i+1}_{N_u})F^{-1}(\bm{u}_{0, N}^{i+1}+\bm{\ell}^{i+1}_{N_u} - D_t\bm{d}_{N_{d_f}}^{i+1}), \bm{v}_{N_u})_{\Omega_f} + \\
+(J\sigma_f(\bm{u}_{0,N}^{i+1}+\bm{\ell}^{i+1}_{N_u}, p_{N_p}^{i+1})F^{-T}, \nabla \bm{v}_{N_u})_{\Omega_f} - (\bm{\lambda}_{N_{\lambda_u}}^{i+1}, \bm{v}_{N_u})_{\Gamma_{FSI}}= (J\bm{b}_f, \bm{v}_{N_u})_{\Omega_f},\\
-(\text{div}(JF^{-1}\bm{u}_{0, N}^{i+1}+\bm{\ell}^{i+1}_{N_u}), q_{N_p})_{\Omega_f} = 0 \\
(\frac{1}{J}\nabla\bm{d}_{N_{d_f}}^{i+1}, \nabla\bm{e}_{N_{d_f}})_{\Omega_f} = (\bm{\lambda}_{N_{\bm{\lambda}_d}}^{i+1}, \bm{e}_{N_{d_f}})_{\Gamma_{FSI}},\\
(\bm{u}_{0, N}^{i+1} +\bm{\ell}^{i+1}_{N_u}, \bm{\mu}_{N_{\lambda_u}})_{\Gamma_{FSI}}-(D_t\bm{d}_{N_{d_s}}^{i+1}, \bm{\mu}_{N_{\lambda_u}})_{\Gamma_{FSI}}=0,\\
(\bm{d}_{N_{d_f}}^{i+1} - \bm{d}_{N_{d_s}}^{i+1}, \bm{\mu}_{N_{\lambda_d}})_{\Gamma_{FSI}}=0.
\end{cases}
\end{equation}

In the previous system, $\bm{\ell}^{i+1}_{N_u}$ is the projection of the finite element lifting function $\bm{\ell}_h^{i+1}$ on the finite dimensional space generated by the velocity reduced basis functions. Once we solve the online system, we can restore the reduced fluid velocity $\bm{u}_{N_u}^{i+1}$ that satisfies the inlet condition $\bm{u}_{N_u}^{i+1}=\bm{u}_{in}(t^{i+1})$ on $\Gamma_{in}$ by using the relation $\bm{u}_{N_u}^{i+1} = \bm{u}_{0, N_u}^{i+1} + \bm{\ell}^{i+1}_{N_u}$.
\section{Results}\label{results mono}
We now present some numerical results that were obtained by adopting a monolithic approach for the toy problem inspired by the Turek--Hron benchmark test case \cite{turek-hron1, turek-hron2}; {in particular, we refer to the test case FSI2 therein, which corresponds to a fluid with Reynolds number $Re=100$.} {We remark again before going any further that, while in the original benchmark problem the Green strain tensor was used for the solid, here, we used the linearized strain tensor, with the reason being that we are not interested in modelling large deformations; however, the values of the parameters are taken from the benchmark FSI2 presented in \cite{turek-hron1, turek-hron2}.}

{All of the numerical simulations for the offline phase were obtained with the use of \textit{multiphenics} \cite{multiphenics}, whereas the online simulations were implemented with \textit{RBniCS} \cite{rbnics}}.

\begin{figure}[H]
\begin{tikzpicture}
\node[anchor=south west,inner sep=0] (image) at (0,0) {\includegraphics[scale=0.4]{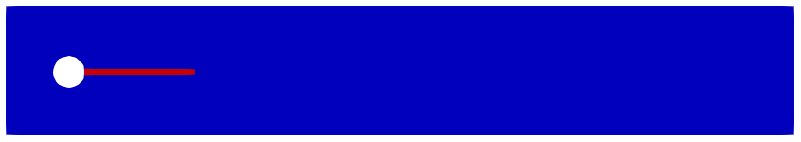}};
\begin{scope}[x={(image.south east)},y={(image.north west)}]
	\draw[white,thick] (0.94,0.5) node {$\Gamma_{out}$};
	\draw[white,thick] (0.13,0.5) node {\small $\Gamma_D^s$};
	\draw[white,thick] (0.03,0.5) node {\small $\Gamma_{in}$};
	\draw[white, thick] (0.5, 1.1) node {$\overset{\mathcal{A}_f(t)}{\curvearrowright}$};
\end{scope}
\end{tikzpicture}
\caption{Reference configuration of the benchmark test case. The solid is depicted in red, while the fluid domain is in blue.\label{turek-hron}}
\end{figure}
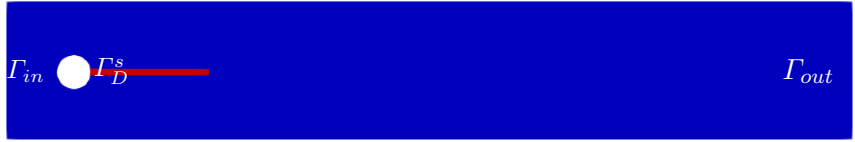

Figure~\ref{turek-hron} represents the physical domain of the problem of interest. The {channel} has a length $L_f=2.5$ cm and a height of $h_f = 0.41$ cm. The cylinder, which is assumed to be at rest and therefore is not considered as part of the solid domain, has center $C=(0.2, 0.2)$ and a radius of $r=0.05$ cm. The deformable bar is $0.35$ cm long and $0.02$ cm thick. 
\vspace{-24pt}

In Table~\ref{values}, we  summarized the values of the physical parameters identifying the fluid and the solid behavior: as we can see, the test case corresponds to a fluid with Reynolds number $Re=100$, {where $Re=\frac{\overline{U}2r}{\nu_f}$, and $r$ is the radius of the immersed cylinder}. We impose homogeneous Dirichlet boundary conditions around the cylinder and on the top and bottom walls of the cavity $\Gamma_{walls}$ for the fluid velocity {and impose no conditions on $\Gamma_N^f$}. We impose a velocity profile at the inlet boundary:
\begin{equation*}
\bm{u}_{in}(t, y):=
\begin{cases}
{\overline{\bm{u}}(0, y)}\frac{1-\text{cos}(\frac{\pi}{2}t)}{2} \quad \text{if }t<2s\\
{\overline{\bm{u}}(0, y)} \quad \text{otherwise},
\end{cases}
\end{equation*}
where
\begin{equation*}
\overline{\bm{u}}(0, y) = 1.5\overline{U}\times \frac{4}{0.1681}y(0.41-y),
\end{equation*}
and the value of $\overline{U}$ is reported in Table~\ref{table monolithic}. 
We also require the bar to be attached to the cylinder; therefore, $\bm{d}_s=0$ on $\Gamma_D^s$.

\begin{table}[H] 
\caption{\centering Values of the physical constants of the fluid and of the solid.\label{values}}
\begin{tabular}{cc}
\toprule
\label{table monolithic}
\textbf{Parameter} &\textbf{Value}\\
\midrule
{Fluid density }$\rho_f$$[10^3\frac{kg}{m^3}]$& $1$\\
{Fluid kinematic viscosity }$\nu_f$$[10^{-3}\frac{m^2}{s}]$& $1$\\
{Mean inflow velocity }$\bar{U}$ $[\frac{m}{s}]$&$1$\\
{Fluid external force }$\bm{b}_f$&$(0, 0)$ \\
{Solid density } $\rho_s$$[10^3\frac{kg}{m^3}]$& $1$\\
{$2$nd Lam\'e constant (solid shear modulus) }$\mu_s$$[10^6\frac{kg}{ms^2}]$& $0.5$\\
{$1$st Lam\'e constant }$\lambda_s$&$0.4$\\
{Solid external force }$\bm{b}_s$ & $(0, 0)$\\
\bottomrule
\end{tabular}
\end{table}

We use a time-step $\Delta t=10^{-2}$ for the discretization in time, and the two constants $\gamma$ and $\beta$ used for the discretization of the structure time derivatives have the following values: $\gamma = 0.25$ and $\beta = 0.5$. The total number of iterations of the simulation is $N_T=10^3$: all of these values are summarized in Table \ref{implementation details}. 
\begin{table}[H] 
\tabcolsep=1.5cm
\caption{Values of the parameters in the time discretization and in the spatial discretization.\label{implementation details}}
\begin{tabular}{cc}
\toprule
\textbf{Time Discretization Parameters} & \textbf{Value}\\
\midrule
Timestep $\Delta T$ & $0.01$s\\
Total number of iterations $N_T$ & $10^3$\\
$\gamma$ & $0.25$\\
$\beta$ & $0.5$\\
\midrule
\textbf{{Space Discretization Parameters} 
} & \textbf{{Value}}\\
\midrule
FE velocity order & $2$\\
FE pressure order & $2$\\
FE displacement order ($\bm{d}_f$ and $\bm{d}_s$) & $2$\\
FE multiplier order ($\bm{\lambda}_u$ and $\bm{\lambda}_d$) & $1$\\
mesh resolution using  \textsf{mshr} mesh generator& $128$\\
\bottomrule
\end{tabular}
\end{table}

Since we are mostly interested in investigating the performance and the ability of the monolithic approach to reproducing the behavior of the coupled system, we adopted a mixed approach: we use a standard FE method, {until the elastic bar starts to oscillate because of the action of the fluid}; then, we run the reduced method. The oscillating behavior of the system takes approximatively $i=800$ iterations to occur: we therefore run the monolithic reduced order method for the remaining $200$ iterations. 

\begin{figure}[H]
\centering
\subfigure[Eigenvalues decay for $u_f$, $s$, $p_f$, $d_f$, and $d_s$.]{\includegraphics[width=0.41\linewidth]{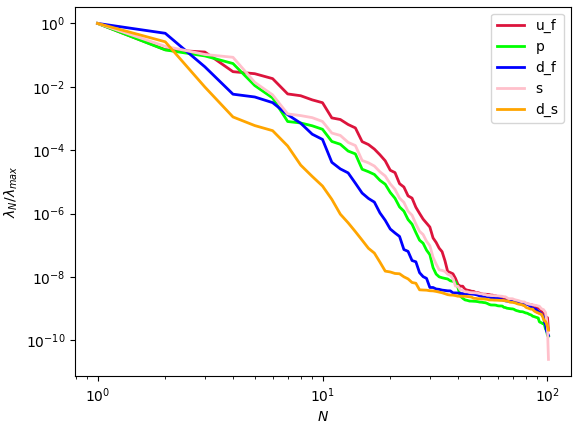}}
\subfigure[Energy retained by the first $100$ modes for $u_f$, $p_f$, $d_f$, and $d_s$.]{\includegraphics[width=0.4\linewidth]{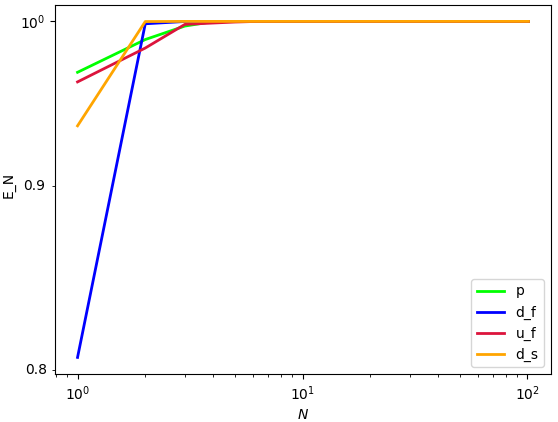}}
\caption{The outcomes of the monolithic POD: eigenvalue decay (\textbf{a}) and retained energy (\textbf{b}) for the first 100 modes.\label{Figure 3}}
\end{figure}

Figure~\ref{Figure 3}a represents the behavior of the first $100$ eigenvalues obtained with the POD on the snapshots of the monolithic system. The eigenvalues of the solid displacement have a faster decay with respect to the others, and this can be justified by the fact that the solid displacement behavior is periodic (the bar oscillates up and down) and that therefore the reduced order model is able to capture this periodic behavior with few modes. On the other hand the fluid velocity eigenvalues present a slower decay with respect to the other components: we expect this to be caused by the fact that, due to the periodic oscillation of the solid, we have the formation of some small vortices in the fluid that propagate into the domain, and this is a more complex phenomenon to reproduce with just a few modes.
In Figure~\ref{Figure 3}b, we can see the behavior of the energy $E_N$ retained by the first $N$ modes for different components of the solution. Here, we give the definition of the retained energy for the fluid velocity component $u_f$, with the understanding that the energy retained for the modes of the other components of the FSI solution is defined in the exact same way:
\begin{equation*}
E_N^u = \frac{\sum_{k=1}^N\lvert\lambda_k^u\rvert}{\sum_{k=1}^{N_u}\lvert\lambda_k^u\rvert}.
\end{equation*}

As we can see from the definition, the energy retained gives us some idea on the amount of information about the physical phenomenon that the first modes carry within: as it has to be expected by looking at the behavior of the eigenvalues, the first modes for the solid retain almost all information about the solid behavior, whereas for the fluid velocity, the energy retained slowly increases to $1$.

\begin{figure}[H]
 
\includegraphics[scale=0.28]{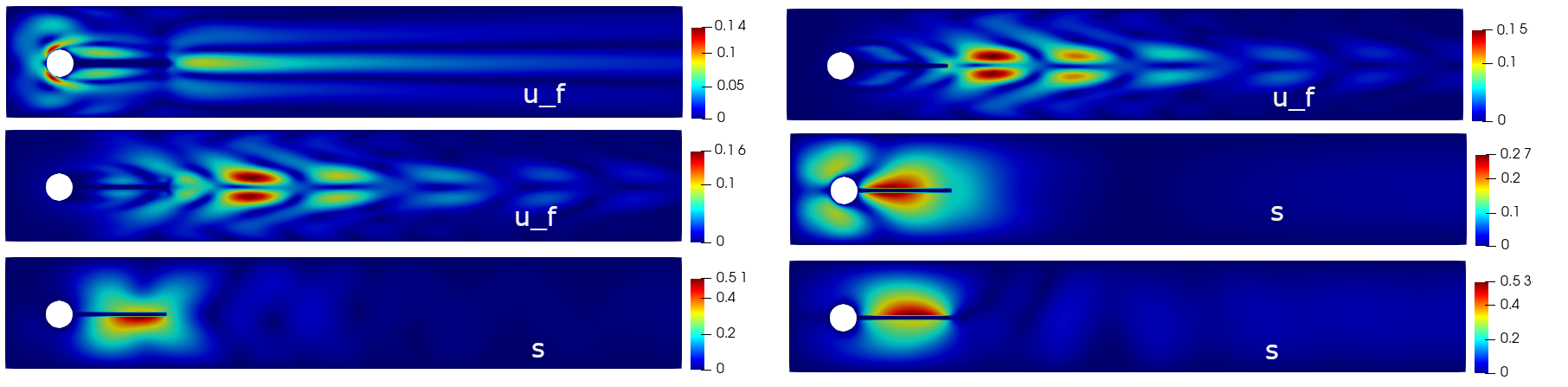}
\caption{\centering The first POD modes for the fluid velocity $u_f$ and for the supremizer $s$ {(magnitude)} .\label{velocity modes monolithic}}
\end{figure}

\begin{figure}[H]
\includegraphics[scale=0.28]{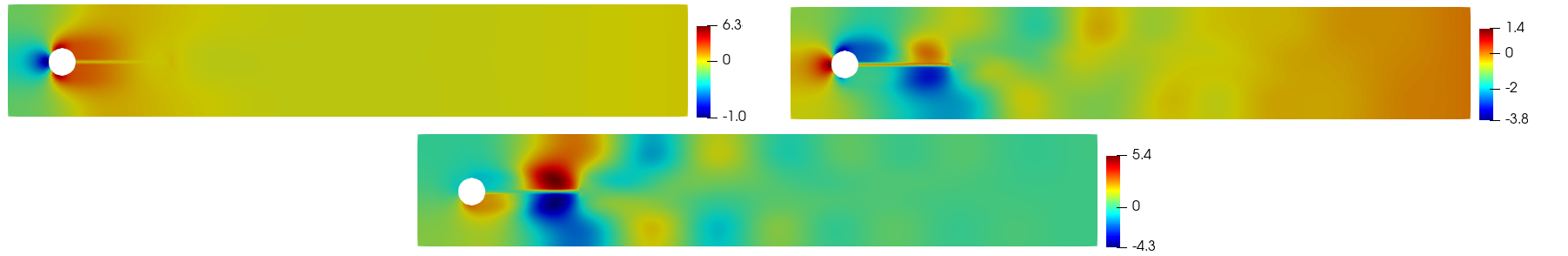}
\caption{The first three POD modes for the fluid {pressure} $p_f$.\label{pressure modes monolithic}}
\end{figure}

\begin{figure}[H]
\includegraphics[scale=0.28]{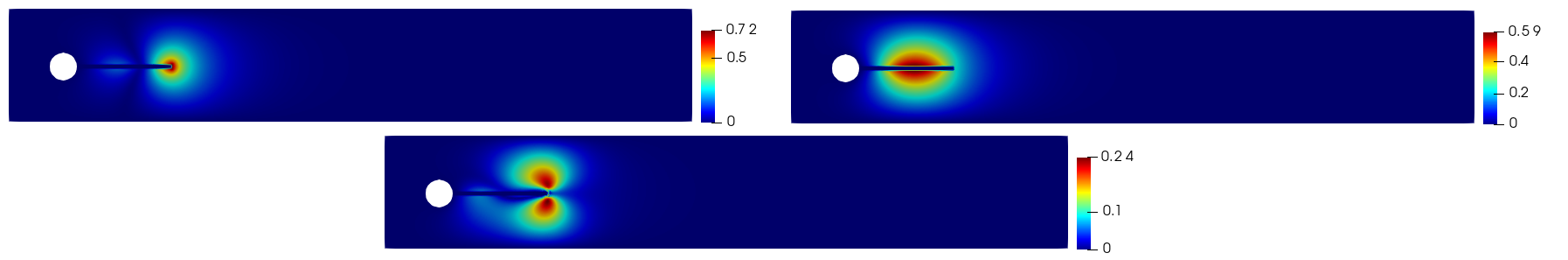}
\caption{The first three POD modes for the fluid displacement $d_f$ {(magnitude)}.\label{mesh modes monolithic}}
\end{figure}

\begin{figure}[H]
\includegraphics[scale=0.28]{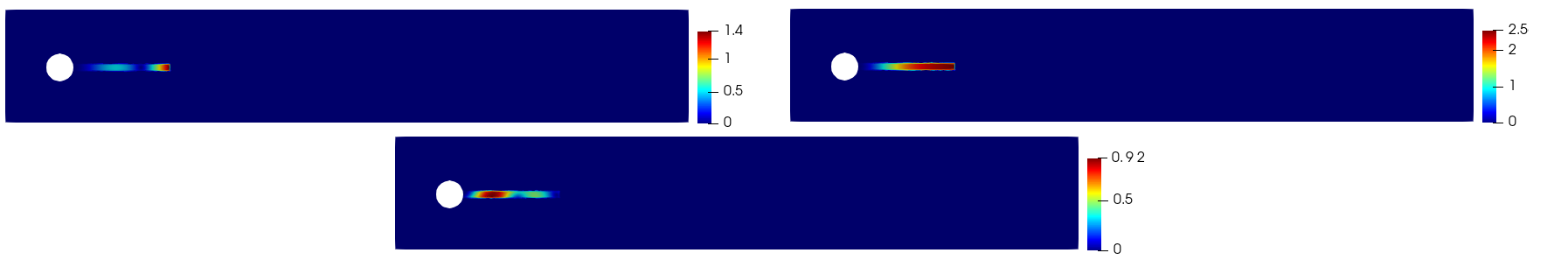}
\caption{The first three POD modes for the solid displacement $d_s$ {(magnitude)}.\label{solid modes monolithic}}
\end{figure}

 Figures~\ref{velocity modes monolithic}--\ref{solid modes monolithic} are intended to help the reader visualize a few outputs of the Proper Orthogonal Decomposition on the snapshots obtained with a monolithic approach. \mbox{Figure~\ref{velocity modes monolithic}} represents the first three modes for the fluid velocity and the first three modes for the supremizer enrichment: we observe how the fluid velocity modes are all zero at the inlet boundary thanks to the lifting function. Figure~\ref{pressure modes monolithic} represents the first three modes for the fluid pressure: as we can see, the modes present a highly oscillating behavior, thus suggesting that the supremizer may be needed in the online phase in order to get rid of any instability in the pressure approximation. Figures~\ref{mesh modes monolithic} and~\ref{solid modes monolithic} represent the modes for the mesh displacement and the solid displacement: we remark that, by looking at the two figures, we can see how the continuity of the displacement along the FSI interface and hence the need for the Lagrange multiplier also in the online phase of the algorithm are not automatically satisfied (as in the partitioned approach on the other hand) by the two sets of modes.

\begin{figure}[H]
\includegraphics[scale=0.4]{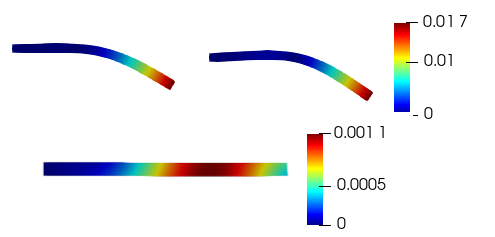}
\caption{The deformation of the bar at time $t=0.9$s s: the FE solution ({top left}) and the reduced order solution ({top right}). Bottom: approximation error $|d_{s, h}-d_{s, N}|$, represented over the solid reference configuration (undeformed state). $N_{d_s}=21$ basis functions were used for the solid displacement. The deformation was magnified by a factor $5$ for visualization {purposes}. \label{monolithic displacement}}
\end{figure}

\begin{figure}[H]
\includegraphics[scale=0.3]{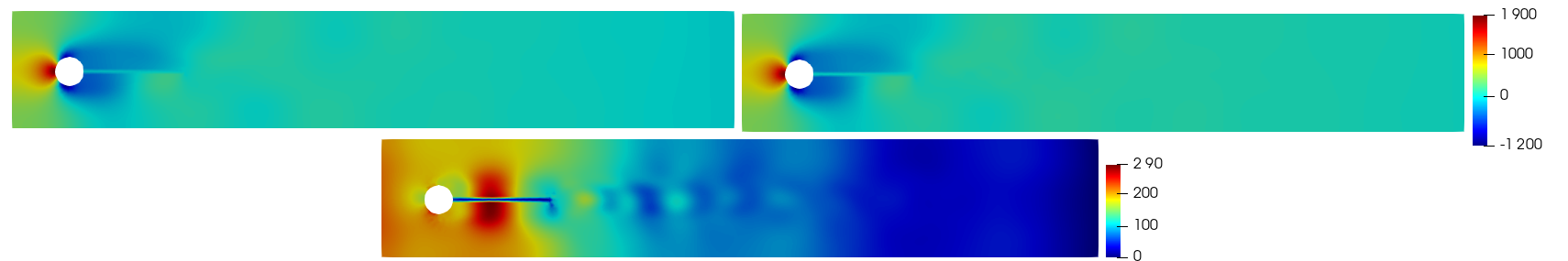}
\caption{Fluid pressure at time $t=0.9$ s: the FE solution ({top left}) and the reduced order solution ({top right}). Bottom: approximation error $|p_{f, h}-p_{f, N}|$. $N_{p}=21$ basis functions were used for the fluid pressure, with the supremizer enrichment {technique}.\label{monolithic pressure}}
\end{figure}

\begin{figure}[H]
\includegraphics[scale=0.3]{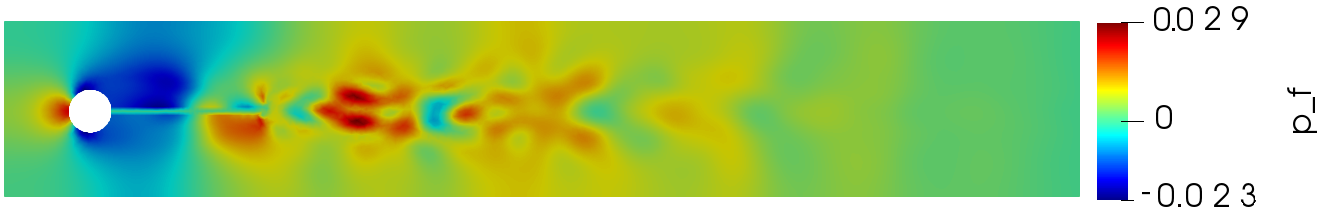}
\caption{Fluid pressure approximation without the implementation of the supremizer enrichment: solution before the {code} 
 {diverges}.\label{pressure no sup}}
\end{figure}

\begin{figure}[H]
\includegraphics[scale=0.3]{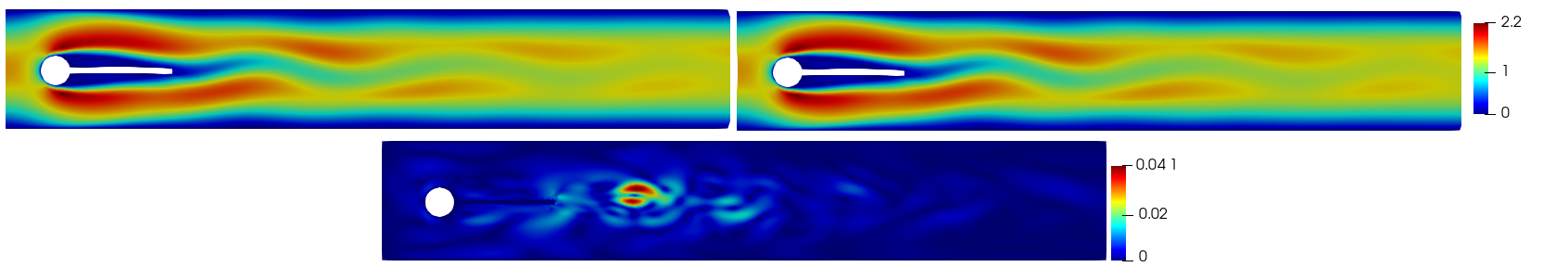}
\caption{Fluid velocity at time $t=0.9$ s: the FE solution ({top left}) and the reduced order solution ({top right}). Bottom: approximation error $|u_{f, h}-u_{f, N}|$. $N_u=42$ basis functions were used for the fluid {velocity}.\label{monolithic velocity}}
\end{figure}

Figure~\ref{monolithic displacement} represents the deformation of the elastic bar at time $t=0.9$ s. The online solution was obtained using $N_{d_s}=21$ basis functions for the component $d_s$ of the solution of the FSI system. As we can see, with this number of reduced bases, the method is able to reconstruct the behavior of the solid component of the coupled system with good accuracy: indeed, the absolute value of the approximation error over the solid domain has a magnitude of $10^{-3}$, which, considered the magnitude of the solid deformation, represents a percentage error of $0.07\%$. In Figures~\ref{monolithic pressure} 
we can see the behavior of the fluid pressure again at time $t=0.9$ s.  The online approximation of the fluid pressure was obtained by employing the supremizer enrichment technique, which requires the introduction of further basis functions in the fluid velocity reduced space: without this technique, the approximation of the fluid pressure becomes unstable and the whole algorithm diverges after a few time--steps, as we can see in Figure~\ref{pressure no sup}. For the fluid pressure, as we can see from Figure~\ref{monolithic pressure} bottom, the approximation error is good, and it represents a percentage error of $0.17\%$.
Figure ~\ref{monolithic velocity} represents the fluid velocity: as we can see, with a Reynolds number of $100$, after some time, we have the developement of some Karman {vortices} that propagate into the fluid domain: the solid bar starts to oscillate and the whole system aquires a periodical behavior. The monolithic algorithm is capable of capturing and reproducing these complex phenomena, such as the Karman {vortices} in the fluid. Figure~\ref{monolithic velocity} shows the behavior of the approximation error: we observe that the error is mostly localized in the region of the domain that is close to the two main vortices that are detached from the solid bar, as it has to be expected, since this represents, from the fluid point of view, the most difficult physical aspect to reproduce.
 
\begin{figure}[H]
\includegraphics[scale=0.6]{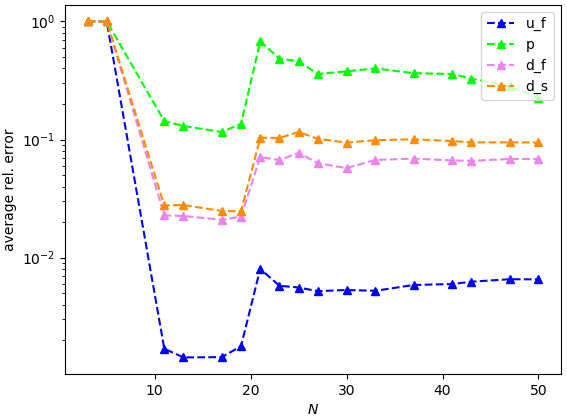}
\caption{{Average relative error as a function of the number of basis functions used in the online~system}.\label{error analysis monolithic}}
\end{figure}

Finally, Figure~\ref{error analysis monolithic} represents the behavior of the average relative error of approximation for the different components of interest of the FSI problem: the average is taken over the number of time-steps. The relative error of approximation is computed in the norm considered in this manuscript for each component of the solution; hence, the $H^1$ seminorm over $\Omega_f$ for $\bm{u}_f$ and $\bm{d}_f$, the $L^2$ norm over $\Omega_f$ for $p_f$, and the $H^1$ seminorm over $\Omega^s$ for $\bm{d}_s$. We would like to remark  that, for this error analysis, we kept the number of reduced basis functions for the Lagrange multipliers fixed, in this specific case, to $N_{\lambda_u}=N_{\lambda_d}=5$. Indeed, we observed an increased presence of stability issues of the online algorithm as we increased the number of modes for the approximation of the Lagrange multipliers. One possible explanation to this is the following: in the online system (see \mbox{Equation~\eqref{online system monolithic}}), the reduced Lagrange multiplier $\bm{\lambda}_{N_{\lambda_u}}$ represents the surface traction created by the \emph{homogenized} fluid velocity $u_{0, N}$ and not by the fluid velocity $u_N$. For this reason, the FE Lagrange multiplier and the reduced order Lagrange multiplier have a different physical interpretation. We suspect this may be the source of instabilities arising by increasing the number of modes for the multipliers: further investigations on this subject will be carried out as future steps.

\begin{figure}[H]
\includegraphics[scale=0.6]{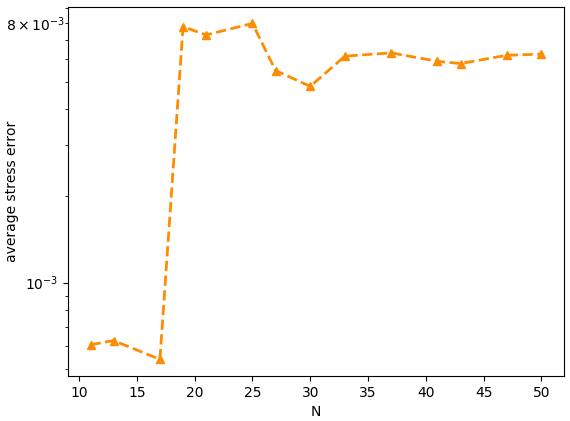}
\caption{Average approximation error for the solid stress at the FSI {interface.} 
\label{average stress error monolithic}}
\end{figure}

It is interesting to see that the average error of approximation decreases up to a certain number of basis function: in this specific case, it decreases up until, more or less, $20$ reduced bases are used. After this threshold, the average relative approximation error increases and reaches a plateau, which indicates the fact that we have no gain in adding more reduced basis functions in the simulation. The increase in the average approximation error indicates that we reached a number of modes after which, if we add more basis functions, we  just add noise to the online system. The same behavior is observed in Figure ~\ref{average stress error monolithic}, which represents the average approximation error of the solid stress at the fluid--structure interface, where the average is taken over the number of time--steps, and the error is computed in the $L^2$ norm: again, the approximation error increases if we use a number of reduced basis that is greater than $20$, confirming the fact that we are just adding noise to the system.
\section{Partitioned Approach}\label{part intro}
In this section, we propose an alternative approach that is instead based on a segregated procedure: the idea is now to solve  the fluid and the solid problems separately and to couple the two physics through some iterative procedure. As we will see, this idea leads to some advantages from the reduction point of view, and  it allows us to work without the employment of additional Lagrange multipliers for the coupling conditions at the fluid--structure interface.
The procedure that we propose is based on a Chorin--Temam projection scheme for the incompressible Navier--Stokes equations \cite{GuermondQuartapelle, Guermond}; in addition, we employ a semi-implicit treatment of the coupling conditions (\cite{BallarinRozzaMaday, Quaini2008, FernandezGerbeauGrandmont}).
\subsection{Time Discretization}\label{part time}
We first present the partitioned scheme after time discretization. $\Delta T$ is the time-step: we also employ an equispaced discretization of the time interval $[0, T]$ in this case. For discretization of the time derivative of a function $f$, we now use first backward difference~BDF1: 
\begin{equation*}
D_tf^{i+1} = \frac{f^{i+1} - f^i}{\Delta T},\qquad D_{tt}f^{i+1}= D_t(D_tf^{i+1}). 
\end{equation*}
The partitioned algorithm reads as follows for $i = 0, \dots, N_T$:
\begin{itemize}
\item Extrapolation of the mesh displacement: 
find $\bm{d}_f^{i+1}\colon\Omega_f\mapsto\mathbb{R}^2$ such that
\begin{equation}
\label{mesh displacement}
\begin{cases}
-\Delta \bm{d}_f^{i+1}=0 \quad\text{in $\Omega_f$,}\\
\bm{d}_f^{i+1} = \bm{d}_s^i \quad\text{on $\Gamma_{FSI}$}. 
\end{cases}
\end{equation}
\item Fluid explicit step:
find $\bm{u}_f^{i+1}\colon\Omega_f\mapsto\mathbb{R}^2$ such that
\begin{equation}
\begin{cases}
J\rho_f \Bigl(\frac{\bm{u}_f^{i+1}-\bm{u}_f^i}{\Delta T} + \nabla \bm{u}_f^{i+1} F^{-1}(\bm{u}_f^{i+1} -D_t \bm{d}_f^{i+1})\Bigr) - \\
-\rho_f\nu_f\text{div}(J\varepsilon(\bm{u}_f^{i+1})F^{-T}) + JF^{-T}\nabla p_f^i= J\bm{b}_f \quad\text{in $\Omega_f$,}\\
\\\bm{u}_f^{i+1} = D_t\bm{d}_f^{i+1} \quad \text{on }\Gamma_{FSI},
\end{cases}
\end{equation}
subject to the Dirichlet boundary condition $\bm{u}_f^{i+1}=\overline{\bm{u}}(t^{i+1})$ on $\Gamma^f_D$, where $\overline{\bm{u}}(t)$ is defined as in Equation \eqref{boundary velocity condition}. {In the above system, $\varepsilon(\bm{u}_f^{i+1})$ is defined as follows:
\begin{equation*}
\varepsilon(\bm{u}_f^{i+1}):=\nabla \bm{u}_f^{i+1}F^{-1} + F^{-T}\nabla^T\bm{u}_f^{i+1}.
\end{equation*}
}
\item Implicit step:
\begin{enumerate}
\item {{Fluid projection substep (pressure Poisson formulation):} 
} find $p_f^{i+1}\colon\Omega_f\mapsto\mathbb{R}^2$ such that
\begin{equation}
\label{pressure implicit step}
\begin{cases}
{-\text{div}(JF^{-1}F^{-T}\nabla p_f^{i+1}) =-\frac{\rho_f}{\Delta t}\text{div}(J F^{-1}\bm{u}_f^{i+1})\quad\text{in $\Omega_f$,}}\\
- F^{-T}\nabla p_f^{i+1}\cdot JF^{-T}\bm{n}_f = \rho_fD_{tt}\bm{d}_s^{i+1}\cdot JF^{-T}\bm{n}_f \quad \text{on $\Gamma_{FSI}$,}
\end{cases}
\end{equation}
subject to the boundary conditions:
\begin{equation}
\label{pressure BC}
p_f^{i+1} = \overline{p}\quad\text{on }\Gamma_{in},
\end{equation}
{where $\overline{p}$ is a prescribed pressure that we computed: a more detailed discussion about this aspect is presented after the structure problem substep.}
\item {{Structure projection substep:}} find $\bm{d}_s^{i+1}\colon\Omega_s\mapsto\mathbb{R}^2$ such that
\begin{equation}
\label{solid implicit step}
\begin{cases}
\rho_sD_{tt}\bm{d}_s^{i+1} - \text{div}P(\bm{d}_s^{i+1})= \bm{b}_s \quad \text{in }\Omega_s,\\
-P(\bm{d}_s^{i+1})\bm{n}_s=J\sigma_f(\bm{u}_f^{i+1}, p_f^{i+1})F^{-T}\bm{n}_f\quad \text{on } \Gamma_{FSI},
\end{cases}
\end{equation}
subject to the boundary condition $\bm{d}_s^{i+1}=0$ on $\Gamma^s_D$.
\end{enumerate}
\end{itemize}

Before going any further, we briefly summarize some important remarks of this formulation of the original FSI problem:
\begin{enumerate}
\item The original Navier--Stokes problem was divided into two subproblems, namely the fluid explicit step and the fluid projection step. In the explicit step, {we take care of the momentum balance of the fluid problem}, whereas in the projection step, we take care of the divergence free condition: this subdivision is a peculiarity of the Chorin--Temam projection scheme; we refer the interested reader to \cite{GuermondQuartapelle, Guermond}. The advantage of adopting such a numerical scheme for the fluid problem is given by the fact that, in this case, we can also use pairs of discrete spaces ($V_h$ and $Q_h$) for the fluid velocity and pressure that do not necessarily satisfy the inf--sup condition; this represents a great advantage for the forthcoming online phase of the method, since we are able to obtain a stable approximation of the fluid pressure also without employing the supremizer enrichment technique, unlike in the monolithic approach.
\item {The treatment of the boundary conditions (in our specific case, the inlet boundary condition) is a delicate aspect of partitioned schemes. If the original problem of interest is provided with a boundary condition for the fluid Cauchy stress tensor $\sigma^f(\bm{u}_f, p_f)\bm{n}_{in}=\bm{g}_{in}$, where $\bm{n}_{in}$ is the normal outgoing the inlet boundary, then during the Chorin--Temam projection scheme, this condition can be \emph{splitted} into two: a natural condition for the fluid velocity explicit step $\varepsilon(\bm{u}_f)\bm{n}_{in}=\bm{0}$ and a Dirichlet condition for the pressure $p_f\bm{n}_{in} = \bm{g}_{in}$ (see for example \cite{Guermond}). However, in our particular toy problem, we have a Dirichlet inlet condition for the velocity: we therefore need a Dirichlet boundary condition also for the pressure in order to obtain uniqueness of solution of the pressure Poisson problem. In this case, we have no specific indication of what value for the pressure to choose at the inlet boundary: for our test case, we decided to compute the inlet pressure value by computing the quantity $\sigma^f(\bm{u}_f, p_f)\bm{n}_{in}$ on the inlet boundary and use the fact that, there, we have $\bm{u}_f=\bm{u}_{in}$.}
\item In the projection step \eqref{pressure implicit step}, we chose a pressure Poisson formulation; it is possible to use a Darcy formulation instead:
find $p_f^{i+1}$ and $\tilde{\bm{u}}_f^{i+1}$ such that
\begin{equation*}
\begin{cases}
\rho_fJ\frac{\tilde{\bm{u}}_f^{i+1}-u_f^{i+1}}{\Delta T} + JF^{-T}\nabla p_f^{i+1} = 0\quad \text{in }\Omega_f,\\
\text{div}(JF^{-1}\tilde{\bm{u}}_f^{i+1}) = 0 \quad \text{in $\Omega_f$}.
\end{cases}
\end{equation*}
However, in view of an efficient model order reduction, we chose to employ a Poisson formulation, since the Darcy formulation requires the introduction of an additional unknown $\tilde{\bm{u}}_f$, which translates in a larger system, comprised of both velocity and pressure, at the implicit step.
\end{enumerate}

In order to enhance the stability of the projection scheme, we employ Robin--Neumann coupling, as proposed in \cite{AstorinoChoulyFernandez, BallarinRozzaMaday}; for other references on this kind of coupling, we refer to \cite{BadiaNobileVergara, FernandezMullaertVidrascu}. We thus replace condition \eqref{pressure implicit step} 
with the following:
\begin{equation}
\label{Robin coupling}
\alpha_{ROB}p^{i+1} + F^{-T}\nabla p^{i+1}\cdot JF^{-T}\bm{n}_f = \alpha_{ROB} p^{i+1,\star} -\rho_fD_{tt}\bm{d}_s^{i+1,\star}\cdot JF^{-T}\bm{n}_f.
\end{equation}

In Equation \eqref{Robin coupling}, $p^{i+1, \star}$ and $\bm{d}_s^{i+1, \star}$ are suitable extrapolations of the fluid pressure and the solid displacement, respectively; we show in the next paragraph which kind of extrapolation we use. The constant $\alpha_{ROB}$ is defined as $\alpha_{ROB} = \frac{\rho_f}{z_p\Delta T}$, where $z_p$ is called the \emph{solid impedance}:
\begin{equation*}
\begin{split}
z_p &= \rho_s c_p,\\
c_p &= \sqrt{\frac{\lambda_s + 2\mu_s}{\rho_s}}.
\end{split}
\end{equation*}
\subsection{Space Discretization}\label{part space}
We now aim to provide the final formulation of the partitioned problem; in order to do so, we need to discretize in space the original problem. We consider hereafter the same function spaces $V^f$, $V^f_0$, $E^f$, and $E^s$ that have been defined in Section~\ref{space discretization section}. {As far as the fluid pressure is concerned, with the Chorin--Temam projection scheme, the solution of the Poisson problem is now in $H^1$, and therefore, we introduce the following pressure function spaces:
\begin{align*}
Q &:= \{q\in H^1(\Omega_f) \text{ st. }q=\overline{p}\text{ on }\Gamma_{in} \},\\
Q_0&:= \{q\in H^1(\Omega_f) \text{ st. }q=0\text{ on }\Gamma_{in} \}.
\end{align*}
}

Additionally, in this case, we discretize the FSI problem in space using second-order Lagrange finite elements for the fluid velocity, the fluid displacement, and the solid displacement, resulting in the discrete spaces $V_h^f\subset V^f$, $V_h^0\subset V^f_0$, $E^f_h\subset E^f$, and $E_h^s\subset E^s$, while the fluid pressure is discretized with first-order Lagrange finite elements, resulting in the discrete spaces $Q_h\subset Q$, $Q_h^0\subset Q_0$ .

We are now ready to present the weak formulation of the original problem for every $i=0, \hdots, N_T$:
\begin{itemize}
\item Extrapolation of the mesh displacement: 
find $\bm{d}_{f, h}^{i+1}\in E_h^f$ such that $\forall \bm{e}_{f,h} \in E^f_h$:
\begin{equation}
\label{mesh displacement}
\begin{cases}
(\nabla\bm{d}_{f, h}^{i+1},\nabla\bm{e}_{f, h})_{\Omega_f} =0\\
\bm{d}_{f, h}^{i+1} = \bm{d}_{s, h}^i \quad\text{on $\Gamma_{FSI}$}. 
\end{cases}
\end{equation}
\item Fluid explicit step:
find $\bm{u}_{f, h}^{i+1}\in V_h^f$ such that $\forall \bm{v}_{f,h}\in V_h^0$:
\begin{equation}
\label{fluid explicit step}
\begin{cases}
\begin{split}
&\rho_f(J\Bigl(\frac{\bm{u}_{f, h}^{i+1}-\bm{u}_{f, h}^i}{\Delta T}\Bigr), \bm{v}_{f,h})_{\Omega_f} 
+ \rho_f(J(\nabla \bm{u}_{f, h}^{i+1}F^{-1}(\bm{u}_{f, h}^{i+1} -D_t\bm{d}_{f, h}^{i+1})),\bm{v}_{f,h})_{\Omega_f}\\
&+ \rho_f\nu_f(J\varepsilon(\bm{u}_{f, h}^{i+1})F^{-T},\nabla \bm{v}_{f,h})_{\Omega_f} + (JF^{-T}\nabla p_{f, h}^i, \bm{v}_{f,h})_{\Omega_f}= (J\bm{b}_f, \bm{v}_{f, h})
\end{split}
\\ \bm{u}_{f, h}^{i+1} = D_t\bm{d}_{f, h}^{i+1} \quad \text{on $\Gamma_{FSI}$,}
\end{cases}
\end{equation}
\item Implicit step: for any $j=0, \dots$ until convergence:
\begin{enumerate}
\item {{Fluid projection substep (pressure Poisson formulation):}} find $p_{f, h}^{i+1,  j+1}\in Q_h$ such that $\forall q_{f,h}\in Q_h^0$:
\begin{myequation}
\begin{split}
&-\frac{\rho_f}{\Delta T}(\text{div}(JF^{-1}\bm{u}_{f, h}^{i+1}), q_{f,h})_{\Omega_f} - \rho_f((D_{tt}\bm{d}_{s, h}^{i+1,  j}), JF^{-T}\bm{n}_fq_{f,h})_{\Gamma_{FSI}} +\\
&+ \alpha_{ROB}(p_{f, h}^{i+1,  j}, q_{f,h})_{\Gamma_{FSI}}= \alpha_{ROB}(p_{f, h}^{i+1,  j+1}, q_{f,h})_{\Gamma_{FSI}} +(JF^{-T}\nabla p_{f, h}^{i+1,  j+1}, F^{-T}\nabla q_{f,h})_{\Omega_f}.
\end{split}
\end{myequation}
\item {{Structure projection substep:}} find $\bm{d}_{s, h}^{i+1,  j+1}\in E^s_h$ such that $\forall \bm{e}_{s, h}\in E^s_h$:
\begin{myequation1}
\rho_s(D_{tt}\bm{d}_{s, h}^{i+1,  j+1}, \bm{e}_{s, h})_{\Omega_f}+ (P(\bm{d}_{s, h}^{i+1,  j+1}), \nabla \bm{e}_{s, h})_{\Omega_s} = -(J\sigma_f(\bm{u}_{f, h}^{i+1}, p_{f, h}^{i+1,  j+1})F^{-T}\bm{n}_f, \bm{e}_{s, h})_{\Gamma_{FSI}} + (\bm{b}_s, \bm{e}_{s, h})_{\Omega_s}
\end{myequation1}
subject to the boundary condition $\bm{d}_{s, h}^{i+1,  j+1} = 0$ on $\Gamma_{D}^s$.
\end{enumerate}
\end{itemize}

We iterate between the two implicit substeps using a fixed point strategy:
\begin{equation}
\label{fixed point criterion}
\text{max}\Bigl(\frac{|| p_{f, h}^{i+1,  j+1} - p_{f, h}^{i+1,  j}||_{Q_h}}{|| p_{f, h}^{i+1,  j+1}||_{Q_h}}; \frac{|| \bm{d}_{s, h}^{i+1,  j+1} - \bm{d}_{s, h}^{i+1,  j}||_{E^s_h}}{||\bm{d}_{s, h}^{i+1,  j+1}||_{E^s_h}}\Bigr) < \varepsilon,
\end{equation}
where $\varepsilon$ is a fixed tolerance.

In the pressure Poisson formulation, to impose the Robin coupling condition, we chose the pressure at the previous implicit iteration, namely $p_f^{i+1,  j}$, as an extrapolation for the fluid pressure, and the same goes for the extrapolation of the structure displacement.
\subsection{Reduced Basis Generation}\label{part pod}
For generation of the reduced basis for the fluid velocity $\bm{u}_{f, h}$ and the fluid displacement $\bm{d}_{f, h}$ we pursue the idea that was first proposed in \cite{BallarinRozzaMaday}. For the solid displacement $\bm{d}_{s, h}$, we employ a standard POD; for fluid pressure $p_{f, h}$, we first introduce a lifting function $\ell_p(t)$ and obtain the homogenized pressure $p_{0, h}(t)=p_{f, h}(t)-\ell_p(t)$ such that $p_{0, h}(t)=0$ on $\Gamma_{in}\times(0,T]$ and then we perform a standard POD. We therefore define the reduced pressure space $Q_N^0:=\text{span}\{\Phi_1^p,\dots,\Phi_{N_p}^p\}$.
\subsubsection{Change of Variable for the Fluid Velocity}
The main idea here is to introduce a change of variable for $\bm{u}_{f, h}$ in the fluid problem in order to transform condition \eqref{fluid explicit step} 
into a homogeneous boundary condition. {The reason for this is 
similar to the reason why we introduced the lifting function in the first place: it is more convinient for the sake of the online system to work with homogeneous boundary conditions. Indeed, even if the reduced basis functions satisfy condition \eqref{fluid explicit step} 
, it is not guaranteed, in general, that an element of the linear space generated by these reduced basis will also satisfy condition \eqref{fluid explicit step} 
. For this reason, we need to introduce a Lagrange multiplier in order to make sure that Equation~\eqref{fluid explicit step}  
is satisfied also at the reduced order level. Therefore, in order to avoid this and in order to design a more efficient reduced method, we chose to transform the non-homogeneous coupling condition into a homogeneous one. First, we transform the non--homogeneous inlet boundary condition $\bm{u}_f^{i+1}={\bm{u}}_{in}(t^{i+1})$ on $\Gamma_{in}$, by introducing a lifting function $\bm{\ell}_u$, similar to that performed for the monolithic approach in Section~\ref{reduced basis monolithic}: we therefore obtain a homogenized fluid velocity $\bm{u}_{0, h}^{i+1}$ such that $\bm{u}_{0, h}^{i+1}=0$ on $\Gamma_{in}$. Afterwards, we define a new variable $\bm{z}_{f, h}^{i+1}$:}\\
\begin{equation*}
\bm{z}_{f, h}^{i+1} := \bm{u}_{0, h}^{i+1} - D_t\bm{d}_{f, h}^{i+1} =  \bm{u}_{f, h}^{i+1} - \bm{\ell}_u^{i+1}- D_t\bm{d}_{f, h}^{i+1}
\end{equation*} 

With this change in variable, Equation \eqref{fluid explicit step} 
is equivalent to the homogeneous boundary condition for the new variable: 
\begin{equation*}
\bm{z}_{f, h}^{i+1}=0 \quad \text{on $\Gamma_{FSI}$},
\end{equation*}
for which no imposition by means of Lagrange multiplier is needed.
Therefore, during the offline phase of the scheme, at every iteration $i+1$, after we compute the homogenized velocity $u_{0, h}^{i+1}$, we compute the change of variable $\bm{z}_{f, h}^{i+1}$. We then consider the following snapshots matrix:
\begin{equation*}
\mathbf{\mathcal{S}}_z = [\bm{z}_{f,h}^1,\hdots,\bm{z}_{f,h}^{N_T}] \in \mathbb{R}^{\mathcal{N}^u_h\times N_T},
\end{equation*}
where $\mathcal{N}^u_h = \text{dim}V_h$. We then apply a POD to the snapshots matrix $\mathbf{\mathcal{S}}_{z}$, and we retain the first $N_{z}$ POD modes $\bm{\Phi}^{z}_1, \hdots, \bm{\Phi}^{z}_{N_{z}}$. We therefore have the reduced space:
\begin{equation*}
V^N:=\text{span}\{\bm{\Phi}^{z}_k\}_{k=1}^{N_{z}},
\end{equation*}
and now it is clear that, since every $\bm{\Phi}^{{z}}_k$ satisfies the condition $\bm{\Phi}^{{z}}_k=0$ on $\Gamma_{FSI}$, then also every element of $V^N$ satisfies the same condition.
\subsubsection{Harmonic Extension of the Fluid Displacement}
In order to generate a reduced basis for the fluid displacement $\bm{d}_f$, we pursue the idea presented in \cite{BallarinRozzaMaday}. Therefore, we start by generating the snapshots matrix related to the solid displacement:
\begin{equation*}
\mathbf{\mathcal{S}}_{d_s} = [\bm{{d}}_{s,h}^1, \hdots, \bm{{{d}}}_{s,h}^{N_T}] \in \mathbb{R}^{\mathcal{N}^{{d}_s}_h\times N_T},
\end{equation*}
where $\mathcal{N}^{{d}_s}_h = \text{dim}E_h^s$, and again, the underline notation denotes the vector of the FE degrees of freedom corresponding to each solution of the solid displacement. We then apply a POD to the snapshots matrix and retain the first $N_{{d}_s}$ POD modes $\bm{\Phi}^{{d}_s}_1, \hdots, \bm{\Phi}^{{d}_s}_{N_{{d}_s}}$, thus defining the reduced space for the solid problem:
\begin{equation*}
E_N^s:=\text{span}\{\bm{\Phi}^{{d}_s}_k\}_{k=1}^{N_{{d}_s}}.
\end{equation*}

We then employ a harmonic extension of each one of the reduced basis $\bm{\Phi}^{{d}_s}_k$ to the fluid domain, thus obtaining the functions $\bm{\Phi}^{{d}_f}_k$ such that
\begin{equation*}
\begin{cases}
-\Delta \bm{\Phi}^{{d}_f}_k = 0 \quad\text{in $\Omega_f$}, \\
\bm{\Phi}^{{d}_f}_k = \bm{\Phi}^{{d}_s}_k \quad\text{on $\Gamma_{FSI}$}.
\end{cases}
\end{equation*}

We can then define the reduced space for the fluid displacement:
\begin{equation*}
E_N^f:= \text{span}\{\bm{\Phi}^{{d}_f}_k\}_{k=1}^{N_{{d}_s}}.
\end{equation*}

The reason for defining the basis functions for $\bm{d}_f$ in such a way instead of employing a standard POD on the set of snapshots for the fluid displacement computed in the offline phase lies in the fact that we want to avoid the introduction of another Lagrange multiplier to impose the non-homogeneous boundary condition \eqref{mesh displacement} 
With our method, we avoid solving the reduced system related to Equation~\eqref{mesh displacement}: indeed, instead of solving an harmonic extension problem at every time-step in the online phase, we solve \emph{once and for all} $N_{{d}_s}$ harmonic extension problems in the expensive offline phase. Then, during the online phase, the reduced fluid displacement is computed just as a linear combination of the basis $\bm{\Phi}^{{d}_f}_k$, with coefficients that are the coefficients of the reduced solid displacement at the previous time-step. We see the final formulation of the online phase of the algorithm in the next section.

\subsection{Online Computational Phase}\label{part online}
We are now ready to present the online formulation of the partitioned procedure. For every $i=0, \hdots, N_T$, we introduce the reduced functions ${z}_{f, N}^{i+1}$, $p_{f, N}^{0, i+1}$, ${d}_{s, N}^{i+1}$ of the following form:
\begin{gather}
\label{z}
{z}_{f, N}^{i+1} = \sum_{k=1}^{N_{{z}_f}}\underline{\textbf{z}}_k^{i+1}\bm{\Phi}^{{z}_f}_k,\\
\label{p}
p_{f, N}^{0, i+1} = \sum_{k=1}^{N_{p}}\underline{p}_k^{0, i+1}\Phi^{p}_k,\\
\label{d}
{d}_{s, N}^{i+1} = \sum_{k=1}^{N_{{d}_s}}\underline{\textbf{d}}_k^{i+1}\bm{\Phi}^{{d}_s}_k.
\end{gather}
Then, the online phase of the partitioned procedure reads as follows:
\subsubsection*{{{Mesh displacement} 
}:}
let $\bm{d}_{f, N}^{i+1}$ be defined by the reduced solid displacement at the previous time-step:
\begin{equation}
\label{d_f}
\bm{d}_{f, N}^{i+1} = \sum_{k=1}^{N_{{d}_s}}\underline{\textbf{d}}_k^i\bm{\Phi}^{{d}_f}_k;
\end{equation}
\subsubsection*{{{Fluid explicit step (with change of variable)}:}}
find $\bm{z}_{f, N}^{i+1}\in V_N$ such that $\forall \bm{v}_{f,N}\in V_N$:
\begin{myequation3}
\begin{split}
&\rho_f\Delta T^{-1}(J({\bm{z}_{f, N}^{i+1}-\bm{u}_{f, N}^i}), \bm{v}_{f,N})_{\Omega_f}
+ \rho_f(J(\nabla(\bm{z}_{f, N}^{i+1}+ D_t\bm{d}_{f, N}^{i+1} + \bm{\ell}_{u, N}^{i+1}){F}^{-1}\bm{z}_{f, N}^{i+1}), \bm{v}_{f,N})_{\Omega_f} + \\
&+ \rho_f\nu_f(J\varepsilon(\bm{z}_{f, N}^{i+1}){F}^{-T},\nabla\bm{v}_{f,N})_{\Omega_f} + (J{F}^{-T}\nabla p_{f, N}^i, \bm{v}_{f,N})_{\Omega_f}
+ \rho_f(J(\nabla\bm{z}_{f, N}^{i+1}{F}^{-1}\bm{\ell}_{f, N}^{i+1}), \bm{v}_{f,N})_{\Omega_f}
= \\
&-\rho_f\Delta T^{-1}(J({D_t\bm{d}_{f, N}^{i+1}}), \bm{v}_{f,N})_{\Omega_f} -\rho_f\Delta T^{-1}(J{\bm{\ell}_{u, N}^{i+1}}, \bm{v}_{f,N})_{\Omega_f}
- \rho_f\nu_f(J\varepsilon(D_t\bm{d}_{f, N}^{i+1}){F}^{-T}, \nabla\bm{v}_{f,N})_{\Omega_f}\\
&- \rho_f\nu_f(J\varepsilon(\bm{\ell}_{f, N}^{i+1}){F}^{-T}, \nabla\bm{v}_{f,N})_{\Omega_f} - \rho_f(J(\nabla(D_t\bm{d}_{f, N}^{i+1} + \bm{\ell}_{u, N}^{i+1}){F}^{-1}\bm{\ell}_{f, N}^{i+1}), \bm{v}_{f,N})_{\Omega_f}
 + (J\bm{b}^f_N, \bm{v}_{f, N})_{\Omega_f},
\end{split}
\end{myequation3}
where again $\bm{b}^f_N$ is the projection of the fluid volume external force on the space $V_N$.
We then restore the reduced fluid velocity: $\bm{u}_{f, N}^{i+1} = \bm{z}_{f, N}^{i+1} + D_t\bm{d}_{f, N}^{i+1} + \bm{\ell}_{u, N}^{i+1}$. Here $\bm{\ell}_{u, N}^{i+1}$ is the projection of the FE lifting function $\ell_{u, h}^{i+1}$ over the reduced basis space $V_N$.
\subsubsection*{{{Implicit step}:}} for any $j=0, \dots$ until convergence:
\begin{enumerate}
\item {{Fluid projection substep:} 
} find $p_{f, N}^{0, i+1, j+1}\in Q_N^0$ such that $\forall q_{f,N}\in Q_N^0$:
\begin{equation*}
\begin{split}
&-\frac{\rho_f}{\Delta T}(\text{div}(JF^{-1}\bm{u}_{f, N}^{i+1}), q_{f,N})_{\Omega_f}- \rho_f((D_{tt}\bm{d}_{s, N}^{i+1, j}), J{F}^{-T}\bm{n}_fq_{f,N})_{\Gamma_{FSI}}\\
&+ \alpha_{ROB}(p_{f, N}^{i+1, j}, q_{f,N})_{\Gamma_{FSI}}- \alpha_{ROB}(\ell_p^{i+1}, q_{f,N})_{\Gamma_{FSI}}\\
&- J{F}^{-T}(\nabla \ell_p^{i+1}, {F}^{-T}\nabla q_{f,N})_{\Omega_f} = \alpha_{ROB}(p_{f, N}^{i+1, j+1}, q_{f,N})_{\Gamma_{FSI}}\\
& +(J{F}^{-T}\nabla p_{f, N}^{i+1, j+1}, {F}^{-T}\nabla q_{f,N})_{\Omega_f};
\end{split}
\end{equation*}
we then recover the reduced fluid pressure $p_{f, N}^{i+1, j+1} = p_{f, N}^{0, i+1, j+1} + \ell_p^{i+1}$.
\item {{Structure projection substep:}} find $\bm{d}_{s, N}^{i+1, j+1}\in E^s_N$ such that $\forall \bm{e}_{s, N}\in E^s_N$:
\begin{myequation2}
\rho_s(D_{tt}\bm{d}_{s, N}^{i+1, j+1},\bm{e}_{s,N})_{\Omega_s} + ({P}(\bm{d}_{s, N}^{i+1, j+1}), \nabla\bm{e}_{s,N})_{\Omega_s} = -(J\sigma_f(\bm{u}_{f, N}^{i+1}, p_{f, N}^{i+1, j+1}){F}^{-T}\bm{n}_f, \bm{e}_{s,N})_{\Gamma_{FSI}} + (\bm{b}^s_N, \bm{e}_{s, N})_{\Omega_s},
\end{myequation2}
where $\bm{b}^s_N$ is the projection of the solid volume external force on the space $E^s_N$.
\end{enumerate}

Additionally, at the reduced order level, the implicit steps are iterated between one another until the stopping criteria Equation~\eqref{fixed point criterion} is satisfied.
\section{Results}\label{part results}
We now present some results obtained by implementing the partitioned algorithm previously described for the same FSI test case, namely the {toy problem inspired by the Turek--Hron benchmark test case: we remark here also that there is a difference with respect to the original test case presented by Turek and Hron, since in \cite{turek-hron1, turek-hron2}, the authors consider the Green strain tensor for the solid, whereas we restored a linear elasticity problem for the structure}. 
The physical constants describing the fluid and the solid properties are the same, and we therefore refer to the values of Table~\ref{values}.

The time-step used for the segregated approach is $\Delta t=10^{-3}$, for a total of $N_T=10^4$ iterations. Additionally, in this case, we adopt a standard FE approach for the first $8000$ iterations and we employ the partitioned reduced order model for the remaining $2000$ iterations. All of the values are reported in Table \ref{parameters partitioned}.
\begin{table}[H] 
\caption{Values of the parameters in the time discretization and in the spatial discretization.\label{parameters partitioned}}
\tabcolsep=1.7cm
\begin{tabular}{cc}
\toprule
\textbf{Time Discretization Parameters} &\textbf{Value}\\
\midrule
$\Delta T$ & $0.0001$ s\\
total number of iterations $N_T$ & $10^4$\\
\midrule
\textbf{{Space Discretization Parameters} 
} &\textbf{{Value}}\\
\midrule
FE velocity order & $2$\\
FE pressure order & $1$\\
FE displacement order & $2$\\
mesh resolution using \textsf{mshr} generator & $128$\\
tolerance $\varepsilon$ for the implicit iterations & $10^{-5}$\\
\bottomrule
\end{tabular}
\end{table}

\begin{figure}[H]
\begin{tabular}{cc}
\includegraphics[scale=0.5]{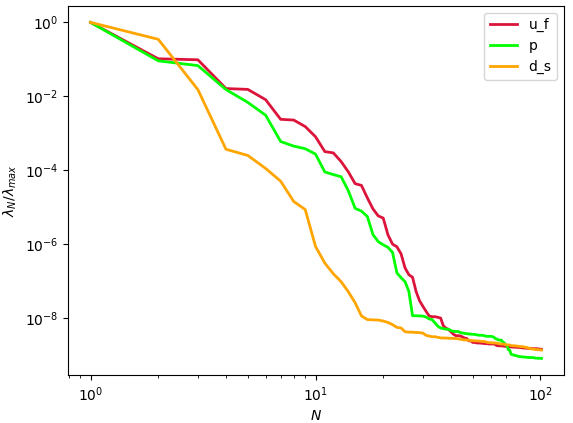}
&\includegraphics[scale=0.5]{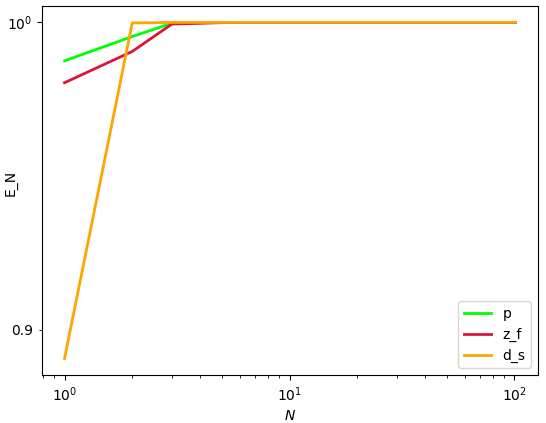}\\
({\bf a}) Decay of the eigenvalues for $\bm{z}_f$, $p_f$, and $d_s$. &({\bf b}) Energy retained for $\bm{z}_f$, $p_f$, and $\bm{d}_s$\\
\end{tabular}
\caption{The outcomes of the partitioned POD: eigenvalue decay (\textbf{a}) and energy retained (\textbf{b}) for the first $100$ modes.}\label{Figure 14}
\end{figure}

{Figure~\ref{Figure 14}a show the rate of decay of the eigenvalues for the fluid pressure $p_f$, the solid displacement $\bm{d}_s$, and the fluid velocity change of variable $\bm{z}_f$. As we can see, in this case, the eigenvalues for the fluid variable $\bm{z}_f$ and for the pressure $p_f$ show almost the same rate of decay, whereas the eigenvalues for the solid displacement decay much faster. In Figure~\ref{Figure 14}b, we can see the energy retained by the first one hundred modes for $\bm{z}_f$, $p_f$ and $\bm{d}_s$. As we can see, the first modes for the fluid pressure are the most energetic ones: they retain almost $10\%$ more energy with respect to the first modes of the solid displacement; additionally, the first modes of the change of variable retain a larger amount of energy with respect to the energy retained by the first modes of the solid displacement.}

\begin{figure}[H]
\includegraphics[scale=0.28]{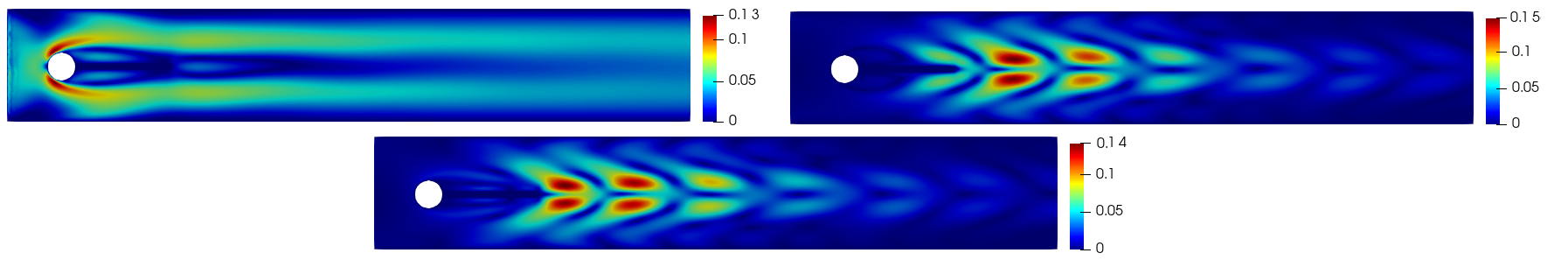}
\caption{{The first three POD modes for the fluid velocity change in variable $z_f$ (magnitude): notice how, for the partitioned approach, the magnitude of the modes is zero not only on the inlet boundary (thanks to the lifting function) but also on the fluid--structure interface (thanks to the implementation of the change of {variable).} 
}\label{modes z_f}}
\end{figure}

\begin{figure}[H]
\includegraphics[scale=0.28]{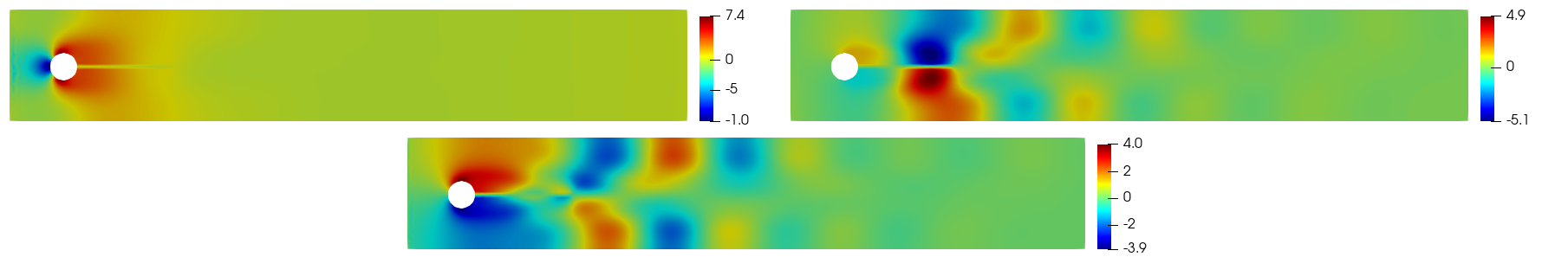}
\caption{{The first three POD modes for the Poisson recovery of the fluid {pressure} $p_f$.}\label{modes p_f partitioned}}
\end{figure}

\begin{figure}[H]
\includegraphics[scale=0.28]{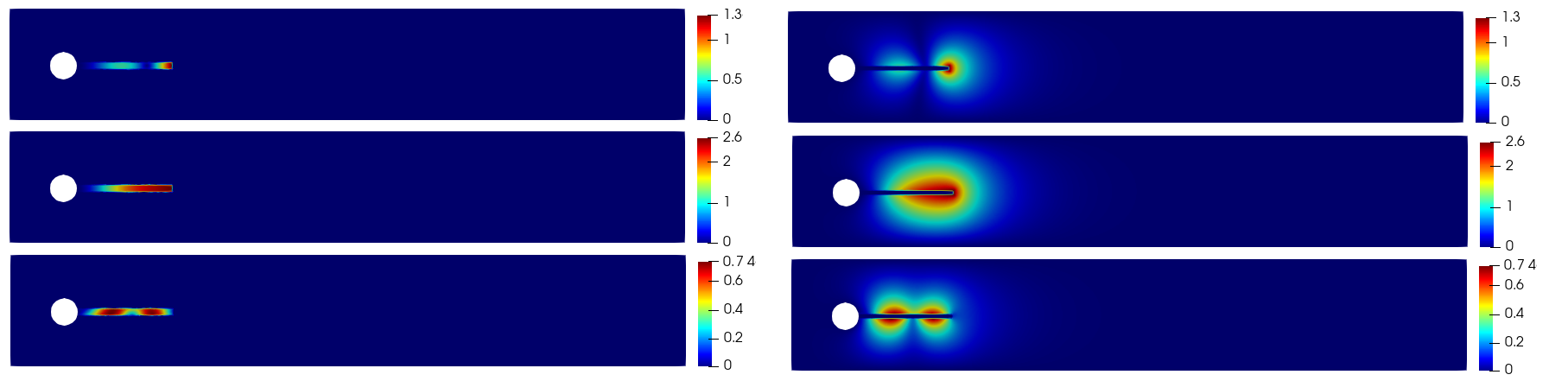}
\caption{{The first three POD modes for the solid displacement $d_s$ ({left column}) and the corresponding mesh displacement modes ({right column}) obtained with an harmonic extension of the basis functions on the left {column.}}\label{modes d_s partitioned}}
\end{figure}

{Figure \ref{modes z_f} represents the first three modes for the fluid change of variable $\bm{z}_f$, and it is interesting to notice how the modes are zero not only on the inlet boundary but also along the FSI interface thanks to the implementation of the change of variable. \mbox{Figure~\ref{modes p_f partitioned}} shows the first three modes for the fluid pressure: as we can see,  in this case, the reduced basis also show a highly oscillatory behavior; nevertheless, thanks to the choice of the Chorin--Temam projection scheme, we are now able to obtain stable approximations of the fluid pressure even without the supremizer enrichment. Finally, Figure \ref{modes d_s partitioned} shows the first three modes for $\bm{d}_s$ (left column); on the right column, we have portraied the corresponding first three basis functions for $\bm{d}_f$, obtained with the harmonic extensions (and not with a POD on the mesh displacement snapshots!).}

\begin{figure}[H]
\includegraphics[scale=0.3]{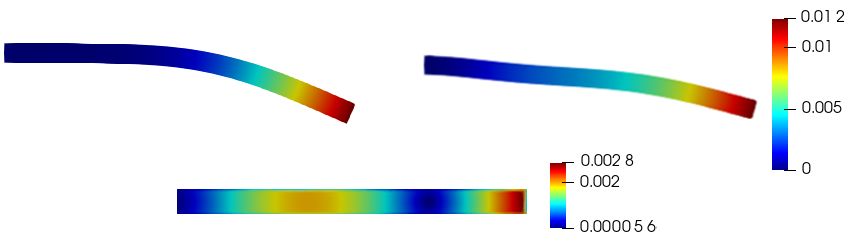}
\caption{Deformation of the structure: FE solution ({top left}), reduced order solution ({top right}), and local approximation error ({bottom}). The approximation was obtained with $N_{d_s}=13$ basis functions. The deformation was magnified by a factor $10$ for visualization {purposes.}\label{displacement partitioned}}
\end{figure}

\begin{figure}[H]
\includegraphics[scale=0.3]{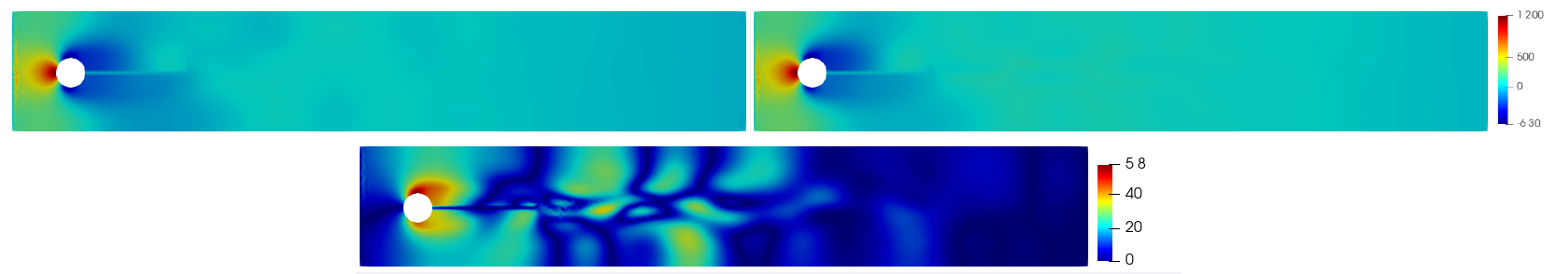}
\caption{Pressure Poisson recovery: FE solution ({top left}), reduced order solution ({top right}), and local approximation error ({bottom}). The approximation was obtained with $N_{p_f}=13$ basis {functions.}\label{pressure partitioned}}
\end{figure}

\begin{figure}[H]
\includegraphics[scale=0.3]{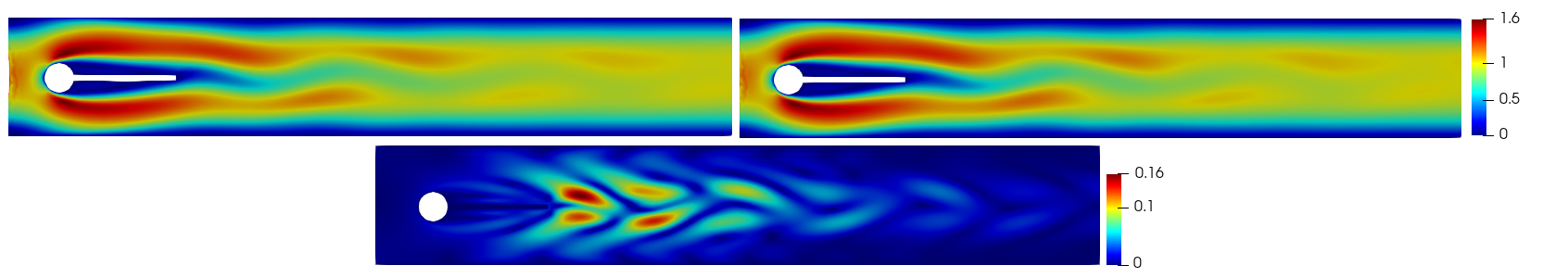}
\caption{Fluid velocity: FE solution ({top left}), reduced order solution ({top right}), and local approximation error ({bottom}). The approximation was obtained with $N_{u_f}=13$ basis {functions.}\label{velocity partitioned}}
\end{figure}

In Figure~\ref{displacement partitioned}, we can see the deformation of the elastic beam at the last time-step of the simulation: {in this case, we used $13$ basis functions for the reduced order approach. We can see that there is some difference between the FE solution and the reduced solution, as the pointwise relative error tells us, with the highest approximation error again at the point where the bar bends downwards: we suspect that this is because, in order to speed up the computations in the online phase, we relaxed the tolerance for the implicit iterations from $10^{-8}$ (FE discretization) to $10^{-5}$. Future investigations will be performed on how to improve the approximation of the solid behavior}. Figure~\ref{pressure partitioned} shows the fluid pressure that was obtained in the partitioned scheme by solving a Poisson problem: this clearly leads to numerical results that are slightly different from the ones in the monolithic approach, but this has to be expected, since we introduced a ``fictious'' Dirichlet boundary condition for the fluid pressure in order to guarantee uniqueness of the solution of the Poisson problem. As we can see, the reduced order fluid pressure accurately represents the FE snapshot even without the use of the supremizer enrichment technique. Finally, in \mbox{Figure~\ref{velocity partitioned}}, we can see the behavior of the fluid velocity: additionally, in this case, we recognize the Karman {vortices} that develop after a while, and in this case, we can see that most of the error is localized in the regions where these vortices detach from the bar and start to propagate into the fluid~domain. 

\begin{figure}[H]
\includegraphics[scale=0.5]{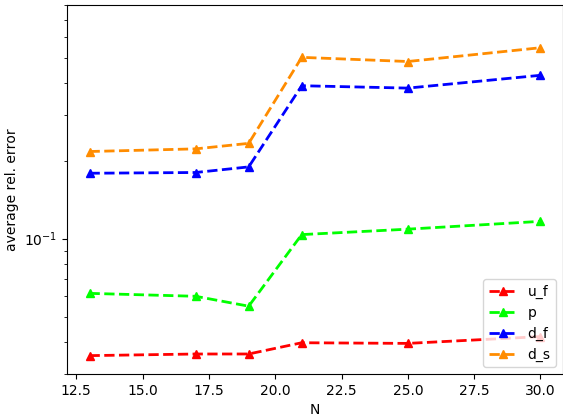}
\caption{Average relative approximation error as a function of the number $N$ of modes used in the online phase. \label{average_approx_error_partitioned}}
\end{figure}

\begin{figure}[H]
\includegraphics[scale=0.5]{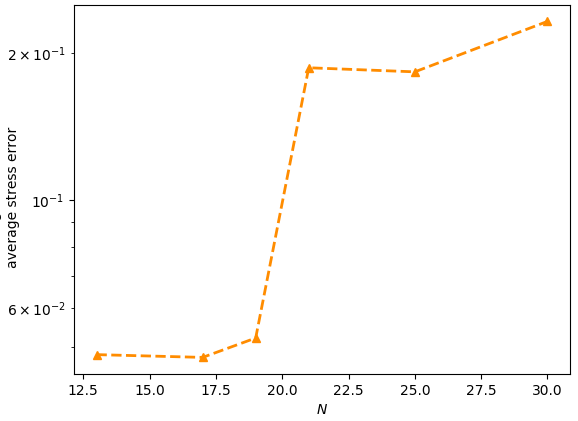}
\caption{Average approximation error for the solid stress at the interface.\label{average stress error partitioned}}
\end{figure}

{In Figure \ref{average_approx_error_partitioned}, we  depicted the behavior of the average relative approximation error (average with respect to time) of the various components of the problem, namely $\bm{u}_f$, $p_f$, $\bm{d}_f$ and $\bm{d}_s$. As we can see, the fluid velocity has the best approximation error compared to the other components, and adding reduced basis functions does not seem to provide a significant improvement in the approximation quality; on the contrary, adding basis functions for $p_f$ does seem to affect the quality of the approximation. As we can see, from $N=19$, it looks like we just add noise to the online system. Finally, we can see that, as expected also from the previous representation of the pointwise approximation error (Figure \ref{displacement partitioned}, the solid displacement is the one for which we have the higher approximation error: again, this could depend on the relaxation of the convergence tolerance of the implicit steps in the online system, but future investigations will be carried out on how to improve the approximation of the solid behavior. Figure \ref{average stress error partitioned} represents the average approximation error of the solid stress at the interface: as we can see, we have a good approximation for the first 19 basis functions, and then again from $N=19$, we get the worst approximation error: again, we suspect this is because we add noise to the system and because, from $N=19$, the average approximation error of all of the components becomes worse.}

\begin{figure}[H]
\includegraphics[scale=0.5]{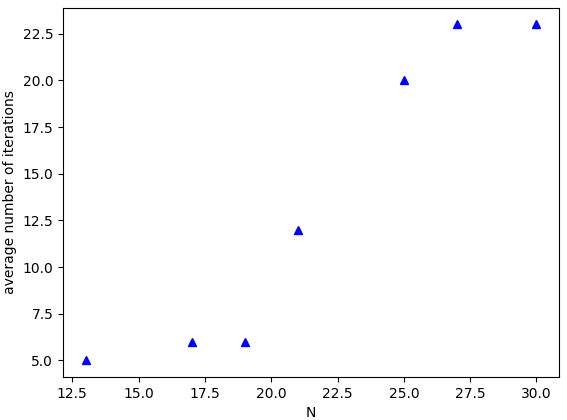}
\caption{Average number of iterations needed in the implicit step in order to reach convergence, as a function of the number $N$ of reduced basis for $\bm{u}_f$, $p_f$, $\bm{d}_s$.\label{number of implicit iterations}}
\end{figure}

{As a final result, we present in Figure \ref{number of implicit iterations} the average number of implicit iterations for the partitioned online phase, plotted against the number of modes used. It is interesting to note that we experience an increase in the average number of iterations in the implicit step (and thus an increase in the computational cost of the online phase) as we increase the number of basis functions; we remember that these results were obtained for a chosen tolerance of $\varepsilon=10^{-5}$ for the fixed point iterations (recall Equation \eqref{fixed point criterion}).}

\section{Discussion}\label{discussion}
{The aim of the work was to provide the reader with an extensive overview about the combination of the RBM with the two main different approaches that are used to address a Fluid--Structure Interaction problem, namely the monolithic and the partitioned approach. We provided insights on the two different kind of algorithms that differ from one another in various aspects: starting from the way the finite element discretization of the original problem is carried out, we then have a significant difference also in the Proper Orthogonal Decomposition. Indeed, even if the underlying idea is to perform a compression of the set of snapshots by means of a POD, we observe that, in the monolithic POD, we work with huge snapshots matrices and with huge inner product matrices because we keep the monolithic structure of the whole approach; in the partitioned approach instead, the snapshots matrices are much smaller. This difference reflects the time that the two PODs take to perform, with the monolithic one being slower than the partitioned one (because we have bigger correlation matrices for which we solve the eigenvalue problem). In addition, the difference between the two approaches also reflects the formulation of the online system: we have a big block system for the monolithic RBM and a series of much smaller systems for the partitioned approach}. 

We  applied the two different model order reduction approaches to a benchmark test case of interest, namely {a test case inspired by} the Turek--Hron FSI test case FSI2 in order to provide the reader with some additional numerical results that provide a better insight on the two procedures.

\begin{figure}[H]
\includegraphics[scale=0.5]{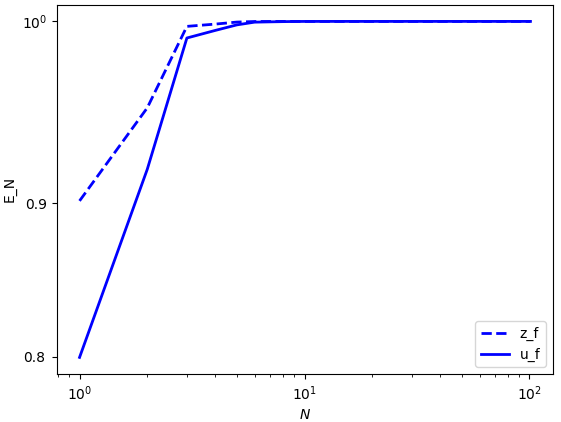}
\caption{A comparison between the energy retained by the first modes for $\bm{z}_f$ and for $\bm{u}_f$\label{retained energy comparison}.}
\end{figure}

Figure~\ref{retained energy comparison} represents a comparison between the amount of energy retained by the first one hundred modes for the fluid velocity $\bm{u}_f$ (modes obtained with the monolithic approach) and the first one hundred modes for the change of variable $\bm{z}_f$ (modes obtained with a partitioned approach). It is interesting to see that the first mode for $\bm{z}_f$ retains almost $10\%$ more energy with respect to the first mode for $\bm{u}_f$, and the general trend is that the first modes for $\bm{z}_f$ retain more energy, thus leading us to believe, at a first glance, that we need fewer reduced basis functions for the fluid momentum equation in the partitioned~approach.

In Table~\ref{final results}, we  summarize some important aspects of both the reduction procedures in order to highlight the differences of the two approaches. {As we can see, the Newton method hwa used in the monolithic approach to solve the huge system and, in the partitioned approach, to solve the fluid explicit step; we used a tolerance for the max. norm of the residual between two consecutive iterations of the solver equal to $6\cdot 10^{-6}$. We remark that it is further possible, if necessary, to speed up the time required for solving the big online system in the monolithic RBM by using some preconditioners: plenty of results on the improvement of the performance of the RBM with preconditioners exist, we refer for example to \cite{DEPARIS2016700}. For the implicit steps of the partitioned online system (pressure Poisson recovery and linear elasticity), we solve the two linear systems with the \textsf{numpy}.\textsf{linalg} solver for linear systems, which is based on a LAPACK routine that performs the LU factorization of the matrix on the left hand side.}
All of the numerical simulations were performed on a computer with $3.50$ GHz per CPU, and the mesh used in the two approaches is the same, as we can see from the mesh resolution reported in \mbox{Tables \ref{implementation details} and \ref{parameters partitioned}}. {For the numerical simulation of the offline phase of both approaches, we relied on \textit{multiphenics} \cite{multiphenics}, which is a Python-based library that helps with the implementation of simulations on conformal meshes of multiphysic problemsM; for the numerical simulation of the online phase, we relied on \textit{rbnics} \cite{rbnics}, which is a Python implementation of several reduced order modelling~techniques.}

\begin{table}
\caption{A comparison between the two approaches.\label{final results}}
\begin{tabular}{cc}
\toprule
\textbf{Monolithic Approach} &\textbf{Value}\\
\midrule
\small $\Delta T$ & $0.01$ s\\
\small number of time iterations $i$ of the RB solver & $200$\\
\small solver for the system & Newton method\\
\small average iterations of Newton method & $4$\\
\small absolute tolerance for Newton method & $||\cdot||_{\infty}<6\cdot 10^{-6}$\\
\small computational time to solve the online system for one time iteration & $224.9$s\\
\midrule
\textbf{{Partitioned Approach} 
} &\textbf{{Value}}\\
\midrule
\small $\Delta T$ & $0.001$ s\\
\small number of time iterations $i$ of the RB solver & $2000$\\
\small solver for explicit fluid step & Newton method\\
\small average number of iterations of Newton method & $2$\\
\small absolute tolerance for Newton method & $||\cdot||_{\infty}<6\cdot 10^{-6}$\\
\small computational time to solve the online systems (explicit + implicit) for one time iteration & $162.68$s\\
\bottomrule
\end{tabular}
\end{table}

The monolithic algorithm, as we have seen, brings along an increase in the number of unknowns to be used in the coupled system: this is because, in order to impose the coupling conditions, we used two Lagrange multipliers, and this leads to an increase in the dimension of the algebraic system to be solved during the online phase. In addition to this, in order to obtain a stable approximation of the reduced order fluid pressure, we adopted a supremizer enrichment of the reduced fluid velocity space: this also leads to an increase of the number of basis functions to be used in the online phase, with a further increase of the dimension of the algebraic system. On the other hand, a monolithic approach is more stable and, therefore, allows for bigger time-steps in the numerical simulations: this turns out to be extremely useful, especially in the case of physical phenomena that take some time to develop, such as the Karman {vortices} in the Turek--Hron benchmark that we considered.
We also make the following remark: for the monolithic approach, one could have used globally continuous spaces for the velocity and the displacement. This approach at the FE level is the one originally adopted by Turek and Hron (see \cite{turek-hron1}), and within the RBM, there are results present in the literature (see \cite{BallarinRozza2016}). At the finite element level, globally continuous spaces for velocity and displacement result in many advantages; however, from the RBM point of view, there are some aspects that we feel are much more easily handled with the approach that we have proposed. First, having globally defined reduced basis functions does not automatically guarantee the balance of the stresses at the interface: this coupling condition would in any case require a weak imposition by means of Lagrange multipliers. A second important aspect is the implementation of the supremizer enrichment technique: if we were to use globally continuous spaces, the supremizer problem would be solved not only in the fluid domain (which is exactly where the supremizer is needed, for stabilization purposes) but also in the entire domain and, hence, in the solid subdomain, adding unnecessary DOFs to the supremizer enrichment problem. Finally, as it is remarked in \cite{BallarinRozza2016}, in the online system, it would not be possible to highlight the contributions corresponding to fluid and solid DOFs, as the unknowns of the online system would be coefficients of a modal expansion (global in both fluid and solid subdomains) and thus not related to a spatial location in the domain. Due to all of these details, we restored to a block formulation from the FE stage: the implementation of the code with this block structure is also much easier thanks to the use of  \textit{multiphenics}.

The partitioned algorithm is very useful, and it gives a lot of control on the systems to be solved during the online phase of the reduction method. In this work we chose a semi-implicit treatment of the coupling conditions: this is, in our opinion, the best choice for the test case considered. Indeed, an implicit treatment of the coupling conditions would be too expensive and would drastically increase the number of sub-iterations at each time-step, both in the offline phase and in the online phase. We also tried to adopt an explicit treatment of the coupling conditions, as suggested for example in \cite{FARHAT20061973}: this approach is unstable for the benchmark considered here, because, due to the slender shape of the domain, the added mass effect plays an important role and leads to an algorithm that diverges after a few time-steps. Therefore, the implicit coupling is the best tradeoff between stability and computational cost. The Chorin--Temam projection scheme has allowed us to work without the need of a supremizer enrichment technique: this gives good control in the number of basis functions to be used in the online phase of the method. In addition to this, the choice of a pressure Poisson formulation allows us to discard the so-called \emph{end of step} velocity and to work just with the intermediate velocity, leading again to a decrease in the dimension of the online system. Finally, the harmonic extension of the fluid displacement basis functions gives the possibility to efficiently compute the mesh displacement in the online phase of the method without the need to solve an additional system. The drawback of a partitioned reduced order model, as we have seen in the Numerical Results section, is that the time-step required in order to have a stable algorithm is, in general, much smaller with respect to the time-step that can be used in a monolithic approach, thus resulting in a larger number of snapshots to be processed in the offline phase. In addition to this, the treatment of the boundary conditions is rather delicate with partitioned approaches and needs to be tailored to the problem at hand, whereas with a monolithic approach, the imposition of these conditions is simpler, either by incorporating them in the weak formulation or by using Lagrange multipliers.

{In this work, we presented numerical results concerning the behavior of the FSI system up to a time of $0.9$ s; however, inspired by the many results present in the literature for the Turek--Hron benchmark test case FSI2, we expect that, given the periodic oscillatory behavior of the structure, for a longer period of observation, vortex shedding phenomena may occur also in our test case. For the Reynolds number considered here, we expect, for both the monolithic and the partitioned approaches, that a much larger number of fluid (pressure and velocity) reduced basis are needed in order to obtain a good approximation of the fluid behavior.  If we were to consider a higher Reynolds number for the fluid, at the reduced order level, the supremizer enirchment technique may not be sufficient to obtain a stable approximation and other stabilization techniques may be required, such as SUPG (see for example \cite{shafqat, Hijazi2020, ali2021reduced}). If the Reynolds number increases significantly and the behavior of the fluid becomes turbulent, then a finite volume based approximation may be the right choice, and in this case, we see an advantage in the partitioned RBM with respect to the monolithic RBM: the partitioned approach indeed potentially allows us to use two different spatial discretization techniques for the fluid and for the solid problem; we refer the interested reader to \cite{Haasdonk, StaHiMoLoRo17, GirfoglioQuainiRozza2019, stabile_stabilized} for some results on the implementation of finite volume discretization within the RBM.

In this work, we also considered a linearized strain tensor for the solid domain because we were interested in developing a procedure for the small deformations range; nonetheless, everything we said can be applied to a nonlinear solid problem as well. If the structure in the system undergoes large deformations, then an ALE formulation may not be the right formalism within which to study the behavior of the coupled system: indeed, it is known that, for large deformation, the ALE formalism may lead to some complications. In this case, we expect a Cut Finite Element approach to be more suited; see for example \cite{PASQUARIELLO2016670}: the investigation of the performance of a Cut Finite Element based RBM for FSI within a monolithic approach is currently under investigation (we refer to \cite{KaNoBaRo, KaBaRO18} for some preliminary results on computational fluid dynamics problems); to the best of our knowledge, there are no results in the literature concerning the application of Cut Finite Element discretization to the RBM within a partitioned approach instead.}

We conclude the discussion with some reflections and perspectives for future works and investigations. For the monolithic reduced order model, even though we introduce a few new variables in the problem formulation (the Lagrange multipliers), we can still make the entire approach computationally feasible: for example, a suggestion could be to perform a rather sparse Proper Orthogonal Decomposition on the first snapshots (when the important physical phenomenon, such as the Karman {vortices} in our test case, is not yet developed) and then to refine the sampling. We are aware that this requires some sort of ``a priori'' knowledge of the physical phenomenon that we are simulating: for the benchmark test case considered, plenty of numerical results already exist that can provide an insight on when to refine the sampling procedure. For other applications, this opens up the possibility of bridging a monolithic reduction procedure with some machine learning algorithm that can be used to investigate, in an efficient way, the overall behavior of the solution to be approximated.

As long as the partitioned approach is concerned, its application is very well indicated for those problems that do not require long time simulations and for industrial applications since the idea of a segregated algorithm is to combine already existing state-of-the-art softwares for computational fluid dynamics and computational solid mechanics. Nonetheless, until alternative stopping criteria and/or alternative treatments of the coupling conditions are further investigated, their application in long time simulations results in an increase in the computational time during the online phase. 
\section{Conclusions}
{In this work, we presented an overview on two possible reduced order models for FSI problems that are based on two different approaches: a monolithic or a partitioned approach. We provided the details of the implementation of the two reduction procedures, both at the FE level and at the reduced order level, analysing the different aspects of the two algorithms, such as the change in variable for the fluid velocity in the partitioned procedure, the creation of mesh displacement basis functions thanks to an harmonic problem, the block structure of the matrices in the monolithic POD, and the different treatment of the coupling conditions at the fluid--structure interface. Finally, we implemented the aforementioned algorithms for a toy problem of interest, which was inspired by the Turek--Hron benchmark test case FSI2. We provided numerical results for the monolithic and for the partitioned RBM, showing, among other things, the behavior of the average approximation error as a function of the number of modes used and the average error between the solid stress at the FSI interface with the FE discretization and with the reduced order discretization. We have seen how the RBM can be modified and adapted in order to be tailored to a monolithic or to a partitioned approach, and the work presented represents an interesting overview, expecially for what concerns segregated reduced order procedures, for which not many results in the literature exist, to the best of our knowledge. 
}
\section*{Acknowledgments}{We acknowledge the support by European Union Funding for Research and Innovation---Horizon 2020 Program---in the framework of European Research Council Executive Agency: Consolidator grant H2020 ERC CoG 2015 AROMA-CFD project 681447 "Advanced Reduced Order Methods with Applications in Computational Fluid Dynamics" (Gianluigi Rozza). We also acknowledge the INDAM-GNCS project "Tecniche Numeriche Avanzate per Applicazioni Industriali". Numerical simulations were obtained with the extension multiphenics of FEniCS \cite{Haasdonk2017,GeeWall} for the high fidelity part and from RBniCS \cite{HesthavenRozzaStamm} for the reduced order part. We acknowledge the developers and contributors of each of the aforementioned libraries. Nonino also acknowledges the support of the Austrian Science Fund (FWF) project F65 ``Taming complexity in Partial Differential Systems'' and the Austrian Science Fund (FWF) project P 33477.}
\bibliographystyle{amsplain}
\bibliography{literature_cutfem.bib}

\end{document}